\documentclass[11pt, a4paper]{article}
\usepackage[cp1251]{inputenc}
\usepackage{amsmath, amscd} \usepackage{euscript}
\usepackage{changepage}
\usepackage{cancel}
\usepackage{tikz-cd}
\usepackage{tikz}
\usepackage{mymatrx5}
\usepackage{mathtools}

\Linestrue\Autonumtrue
\usepackage{mathrsfs}
\usepackage{amssymb}
\usepackage[matrix, arrow, curve]{xy}
\begin{document}
	
	\newcounter{bnomer} \newcounter{snomer}
	\newcounter{bsnomer}
	\setcounter{bnomer}{0}
	\renewcommand{\thesnomer}{\thebnomer.\arabic{snomer}}
	\renewcommand{\thebsnomer}{\thebnomer.\arabic{bsnomer}}
	\renewcommand{\refname}{\begin{center}\large{\textbf{References}}\end{center}}
	
	\setcounter{MaxMatrixCols}{14}
	
	\newcommand\restr[2]{{
			\left.\kern-\nulldelimiterspace 
			#1 
			\right|_{#2} 
	}}
	
	\newcommand{\sect}[1]{%
		\setcounter{snomer}{0}\setcounter{bsnomer}{0}
		\refstepcounter{bnomer}
		\par\bigskip\begin{center}\large{\textbf{\arabic{bnomer}. {#1}}}\end{center}}
	\newcommand{\sectappendix}[1]{%
		\setcounter{snomer}{0}
		\par\bigskip\begin{center}\large{\textbf{A. {#1}}}\end{center}}
	\newcommand{\sst}[1]{%
		\refstepcounter{bsnomer}
		\par\bigskip\textbf{\arabic{bnomer}.\arabic{bsnomer}. {#1}.}}
	\newcommand{\sstappendix}[1]{%
		\refstepcounter{bsnomer}
		\par\bigskip\textbf{A.\arabic{bsnomer}. {#1}.}}
	\newcommand{\lemmpappendix}[2]{%
		\refstepcounter{snomer}
		\par\medskip\textbf{Lemma A.\arabic{snomer}. }{\emph{#1}}
		\par\textsc{Proof}. {#2}\hspace{\fill}$\square$\par\medskip}
	\newcommand{\theopappendix}[2]{%
		\refstepcounter{snomer}
		\par\textbf{Theorem A.\arabic{snomer}. }{\emph{#1}}
		\par\textsc{Proof}. {#2}\hspace{\fill}$\square$\par}
	\newcommand{\proppappendix}[2]{%
		\refstepcounter{snomer}
		\par\medskip\textbf{Proposition A.\arabic{snomer}. }{\emph{#1}}
		\par\textsc{Proof}. {#2}\hspace{\fill}$\square$\par\medskip}
	\newcommand{\defiappendix}[1]{%
		\refstepcounter{snomer}
		\par\medskip\textbf{Definition A.\arabic{snomer}. }{#1}\par\medskip}
	\newcommand{\notaappendix}[1]{%
		\refstepcounter{snomer}
		\par\medskip\textbf{Remark A.\arabic{snomer}. }{#1}\hspace{\fill}$\bigcirc$\par\medskip}
	
	\newcommand{\examappendix}[2]{%
		\refstepcounter{snomer}
		\par\textbf{Example A.\arabic{snomer}. }{#2} {{#1}}\hspace{\fill}$\bigcirc$\par}
	
	\newcommand{\defi}[1]{%
		\refstepcounter{snomer}
		\par\medskip\textbf{Definition \arabic{bnomer}.\arabic{snomer}. }{#1}\par\medskip}
	\newcommand{\theo}[2]{%
		\refstepcounter{snomer}
		\par\textbf{Theorem \arabic{bnomer}.\arabic{snomer}. }{#2} {\emph{#1}}\hspace{\fill}$\square$\par}
	\newcommand{\mtheop}[2]{%
		\refstepcounter{snomer}
		\par\textbf{Theorem \arabic{bnomer}.\arabic{snomer}. }{\emph{#1}}
		\par\textsc{Proof}. {#2}\hspace{\fill}$\square$\par}
	\newcommand{\mcorop}[2]{%
		\refstepcounter{snomer}
		\par\textbf{Corollary \arabic{bnomer}.\arabic{snomer}. }{\emph{#1}}
		\par\textsc{Proof}. {#2}\hspace{\fill}$\square$\par}
	\newcommand{\mtheo}[1]{%
		\refstepcounter{snomer}
		\par\medskip\textbf{Theorem \arabic{bnomer}.\arabic{snomer}. }{\emph{#1}}\par\medskip}
	\newcommand{\theoc}[2]{%
		\refstepcounter{snomer}
		\par\medskip\textbf{Theorem \arabic{bnomer}.\arabic{snomer}. }{#1} {\emph{#2}}\par\medskip}
	\newcommand{\mlemm}[1]{%
		\refstepcounter{snomer}
		\par\medskip\textbf{Lemma \arabic{bnomer}.\arabic{snomer}. }{\emph{#1}}\par\medskip}
	\newcommand{\mprop}[1]{%
		\refstepcounter{snomer}
		\par\medskip\textbf{Proposition \arabic{bnomer}.\arabic{snomer}. }{\emph{#1}}\par\medskip}
	\newcommand{\theobp}[2]{%
		\refstepcounter{snomer}
		\par\textbf{Theorem \arabic{bnomer}.\arabic{snomer}. }{#2} {\emph{#1}}\par}
	\newcommand{\theop}[2]{%
		\refstepcounter{snomer}
		\par\textbf{Theorem \arabic{bnomer}.\arabic{snomer}. }{\emph{#1}}
		\par\textsc{Proof}. {#2}\hspace{\fill}$\square$\par}
	\newcommand{\theopwop}[2]{%
		\refstepcounter{snomer}
		\par\textbf{Theorem \arabic{bnomer}.\arabic{snomer}. }{\emph{#1}}}
	\newcommand{\theosp}[2]{%
		\refstepcounter{snomer}
		\par\textbf{Theorem \arabic{bnomer}.\arabic{snomer}. }{\emph{#1}}
		\par\textbf{Sketch of the proof}. {#2}\hspace{\fill}$\square$\par}
	\newcommand{\exam}[1]{%
		\refstepcounter{snomer}
		\par\medskip\textbf{Example \arabic{bnomer}.\arabic{snomer}. }{#1}\par\medskip}
	\newcommand{\deno}[1]{%
		\refstepcounter{snomer}
		\par\textbf{Notation \arabic{bnomer}.\arabic{snomer}. }{#1}\par}
	\newcommand{\lemm}[1]{%
		\refstepcounter{snomer}
		\par\textbf{Lemma \arabic{bnomer}.\arabic{snomer}. }{\emph{#1}}\hspace{\fill}$\square$\par}
	\newcommand{\lemmp}[2]{%
		\refstepcounter{snomer}
		\par\medskip\textbf{Lemma \arabic{bnomer}.\arabic{snomer}. }{\emph{#1}}
		\par\textsc{Proof}. {#2}\hspace{\fill}$\square$\par\medskip}
	\newcommand{\coro}[1]{%
		\refstepcounter{snomer}
		\par\textbf{Corollary \arabic{bnomer}.\arabic{snomer}. }{\emph{#1}}\hspace{\fill}$\square$\par}
	
	\newcommand{\coropappendix}[2]{%
		\refstepcounter{snomer}
		\par\medskip\textbf{Corollary A.\arabic{snomer}. }{\emph{#1}}
		\par\textsc{Proof}. {#2}\hspace{\fill}$\square$\par\medskip}
	
	\newcommand{\corowop}[1]{%
		\refstepcounter{snomer}
		\par\textbf{Corollary \arabic{bnomer}.\arabic{snomer}. }{\emph{#1}}}
	\newcommand{\mcoro}[1]{%
		\refstepcounter{snomer}
		\par\textbf{Corollary \arabic{bnomer}.\arabic{snomer}. }{\emph{#1}}\par\medskip}
	\newcommand{\corop}[2]{%
		\refstepcounter{snomer}
		\par\textbf{Corollary \arabic{bnomer}.\arabic{snomer}. }{\emph{#1}}
		\par\textsc{Proof}. {#2}\hspace{\fill}$\square$\par}
	\newcommand{\nota}[1]{%
		\refstepcounter{snomer}
		\par\medskip\textbf{Remark \arabic{bnomer}.\arabic{snomer}. }{#1}\hspace{\fill}$\bigcirc$\par\medskip}
	\newcommand{\propp}[2]{%
		\refstepcounter{snomer}
		\par\medskip\textbf{Proposition \arabic{bnomer}.\arabic{snomer}. }{\emph{#1}}
		\par\textsc{Proof}. {#2}\hspace{\fill}$\square$\par\medskip}
	\newcommand{\hypo}[1]{%
		\refstepcounter{snomer}
		\par\medskip\textbf{Conjecture \arabic{bnomer}.\arabic{snomer}. }{\emph{#1}}\par\medskip}
	\newcommand{\prop}[1]{%
		\refstepcounter{snomer}
		\par\textbf{Proposition \arabic{bnomer}.\arabic{snomer}. }{\emph{#1}}\hspace{\fill}$\square$\par}
	
	\newcommand{\proof}[2]{%
		\par\medskip\textsc{Proof{#1}}. {#2}\hspace{\fill}$\square$\par\medskip}
	
	\newcommand{\okr}[2]{%
		\refstepcounter{snomer}
		\par\medskip\textbf{{#1} \arabic{bnomer}.\arabic{snomer}. }{\emph{#2}}\par\medskip}
	
	\newcommand{\hooklongrightarrow}{\lhook\joinrel\longrightarrow}

	\newcommand{\Ind}[3]{%
		\mathrm{Ind}_{#1}^{#2}{#3}}
	\newcommand{\Res}[3]{%
		\mathrm{Res}_{#1}^{#2}{#3}}
	\newcommand{\epsi}{\epsilon}
	\newcommand{\tri}{\triangleleft}
	\newcommand{\Supp}[1]{%
		\mathrm{Supp}(#1)}
	
	\newcommand{\lee}{\leqslant}
	\newcommand{\gee}{\geqslant}
	\newcommand{\reg}{\mathrm{reg}}
	\newcommand{\Ann}{\mathrm{Ann}\,}
	\newcommand{\Cent}[1]{\mathbin\mathrm{Cent}({#1})}
	\newcommand{\PCent}[1]{\mathbin\mathrm{PCent}({#1})}
	\newcommand{\Exp}[1]{\mathbin\mathrm{Exp}({#1})}
	\newcommand{\empr}[2]{[-{#1},{#1}]\times[-{#2},{#2}]}
	\newcommand{\sreg}{\mathrm{sreg}}
	\newcommand{\ilm}{\varinjlim}
	\newcommand{\plm}{\varprojlim}
	\newcommand{\codim}{\mathrm{codim}\,}
	\newcommand{\chara}{\mathrm{char}\,}
	\newcommand{\rk}{\mathrm{rk}\,}
	\newcommand{\chr}{\mathrm{ch}\,}
	\newcommand{\Ker}{\mathrm{Ker}\,}
	\newcommand{\Pic}{\mathrm{Pic}\,}
	\newcommand{\Imm}{\mathrm{im}\,}
	\newcommand{\id}{\mathrm{id}}
	\newcommand{\Ad}{\mathrm{Ad}}
	\newcommand{\col}{\mathrm{col}}
	\newcommand{\row}{\mathrm{row}}
	\newcommand{\low}{\mathrm{low}}
	\newcommand{\pho}{\hphantom{\quad}\vphantom{\mid}}
	\newcommand{\fho}[1]{\vphantom{\mid}\setbox0\hbox{00}\hbox to \wd0{\hss\ensuremath{#1}\hss}}
	\newcommand{\wt}{\widetilde}
	\newcommand{\wh}{\widehat}
	\newcommand{\ad}[1]{\mathrm{ad}_{#1}}
	\newcommand{\tr}{\mathrm{tr}\,}
	\newcommand{\GL}{GL}
	\newcommand{\PGL}{PGL}
	\newcommand{\Pro}{P}
	\newcommand{\SO}{SO}
	\newcommand{\Orth}{O}
	\newcommand{\Sp}{Sp}
	\newcommand{\Gra}{Gr}
	\newcommand{\GraO}{GrO}
	\newcommand{\GraS}{GrS}
	\newcommand{\Fl}{Fl}
	\newcommand{\FlO}{FlO}
	\newcommand{\FlS}{FlS}
	\newcommand{\St}{St}

	\newcommand{\SL}{SL}

	\newcommand{\Sa}{\mathrm{S}}
	\newcommand{\Ua}{\mathrm{U}}
	\newcommand{\Mat}{\mathrm{Mat}}
	\newcommand{\Stab}{\mathrm{Stab}}
	\newcommand{\htt}{\mathfrak{h}}
	\newcommand{\spt}{\mathfrak{sp}}
	\newcommand{\slt}{\mathfrak{sl}}
	\newcommand{\sot}{\mathfrak{so}}

	\newcommand{\vfi}{\varphi}
	\newcommand{\teta}{\vartheta}
	\newcommand{\Bfi}{\Phi}
	\newcommand{\Fp}{\mathbb{F}}
	\newcommand{\Rp}{\mathbb{R}}
	\newcommand{\Zp}{\mathbb{Z}}
	\newcommand{\Cp}{\mathbb{C}}
	\newcommand{\Np}{\mathbb{N}}
	\newcommand{\ut}{\mathfrak{u}}
	\newcommand{\at}{\mathfrak{a}}
	\newcommand{\hei}{\mathfrak{hei}}
	\newcommand{\nt}{\mathfrak{n}}
	\newcommand{\mt}{\mathfrak{m}}
	\newcommand{\rt}{\mathfrak{r}}
	\newcommand{\rad}{\mathfrak{rad}}
	\newcommand{\bt}{\mathfrak{b}}
	\newcommand{\gt}{\mathfrak{g}}
	\newcommand{\vt}{\mathfrak{v}}
	\newcommand{\pt}{\mathfrak{p}}
	\newcommand{\Xt}{\mathfrak{X}}
	\newcommand{\Po}{\mathcal{P}}
	\newcommand{\PV}{\mathcal{PV}}
	\newcommand{\Uo}{\EuScript{U}}
	\newcommand{\Fo}{\EuScript{F}}
	\newcommand{\Do}{\EuScript{D}}
	\newcommand{\Eo}{\EuScript{E}}
	\newcommand{\Xu}{\mathcal{X}}
	\newcommand{\Iu}{\mathcal{I}}
	\newcommand{\Mo}{\mathcal{M}}
	\newcommand{\Nu}{\mathcal{N}}
	\newcommand{\Gu}{\mathcal{G}}
	\newcommand{\Ro}{\mathcal{R}}
	\newcommand{\Co}{\mathcal{C}}
	\newcommand{\Lo}{\mathcal{L}}
	\newcommand{\Ou}{\mathcal{O}}
	\newcommand{\Uu}{\mathcal{U}}
	\newcommand{\Au}{\mathcal{A}}
	\newcommand{\Vu}{\mathcal{V}}
	\newcommand{\Du}{\mathcal{D}}
	\newcommand{\Bu}{\mathcal{B}}
	\newcommand{\Sy}{\mathcal{Z}}
	\newcommand{\Sb}{\mathcal{F}}
	\newcommand{\Gr}{\mathcal{G}}
	\newcommand{\rtc}[1]{C_{#1}^{\mathrm{red}}}

	\newcommand{\JSpec}[1]{\mathrm{JSpec}\,{#1}}
	\newcommand{\PSpec}[1]{\mathrm{PSpec}\,{#1}}
	\newcommand{\APbr}[1]{\mathrm{span}\{#1\}}
	\newcommand{\APbre}[1]{\langle #1\rangle}
	\newcommand{\APro}[1]{\setcounter{AP}{#1}\Roman{AP}}\newcommand{\apro}[1]{{\rm\setcounter{AP}{#1}\roman{AP}}}
	\newcommand{\ot}{\xleftarrow[]{}}
	\newcounter{AP}

	\author{Mikhail Ignatev and Ivan Penkov}
	\date{}
	\title{Automorphism groups of ind-varieties of generalized flags}\maketitle
	\begin{abstract} 
		We compute the group of automorphisms of an arbitrary ind-variety of (possibly isotropic) generalized flags. Such an ind-variety is a homogeneous ind-space for one of the ind-groups $SL(\infty)$, $O(\infty)$ or $Sp(\infty)$. We show that the respective automorphism groups are much larger than $SL(\infty)$, $O(\infty)$ or $Sp(\infty)$, and present the answer in terms of Mackey groups. The latter are groups of automorphisms of nondegenerate  pairings of (in general infinite-dimensional) vector spaces. An explicit matrix form of the automorphism group of an arbitrary ind-variety of generalized flags  is also given. The case of the Sato grassmannian is considered in detail, and its automorphism group is the projectivization of the connected component of unity in the group known as Japanese $\GL(\infty)$.   
		
		
	\end{abstract}
	
	\begin{center}
		{\large \bfseries Introduciton}
	\end{center}
	
	Given a homogenous space $X$, it is a natural problem to compute its automorphism group $\mathrm{Aut}\,X$. In the case when $X$ is a complex flag variety, that is, $X=G/P$ for a connected reductive complex algebraic group and a parabolic subgroup $P\subset G$, the automorphism group of $X$ is well known. Moreover, it is a classical result that here the connected component $\mathrm{Aut}^{0}\, X$ of the identity equals the projectivized group $\Pro G$, except in some special cases as described in \cite{Onishchik}.
	
	In this paper, we would like to pose and solve the problem of computing $\mathrm{Aut}\, X$ for a class of homogeneous ind-varieties $X$. This is the class of ind-varieties of generalized flags introduced by that name in \cite{DimitrovPenkov} but also considered earlier in several works, see for instance \cite{DmitrovPenkovWolf}, \cite{NatarajanRodriguez-CarringtonWolf}. These ind-varieties can be defined simply as $G/P$ where $G$ is one of the ind-groups $\SL(\infty)=\ilm \SL(n)$, $\SO(\infty)=\ilm \SO(n)$, $\Sp(\infty)=\ilm \Sp(2n)$ and $P$ is a splitting parabolic subgroup, i.e., a subgroup for which the intersections $P\cap \SL(n)$, $P\cap \SO(2n)$, $P\cap \SO(2n+1)$, $P\cap \Sp(2n)$ are parabolic subgroups of $\SL(n)$, $\SO(2n)$, $\SO(2n+1)$, $\Sp(2n)$ for all $n$, respectively. 
	The definition from \cite{DimitrovPenkov} can be considered as a flag realization of the ind-varieties $G/P$ as above, and is recalled in Section \ref{section:brief_background} below. The main idea of that approach is that one designates certain chains of subspaces in the natural representation $V$ of $\SL(\infty)$ as generalized flags, and then defines an ind-variety of generalized flags as the ind-variety of generalized flags which differ only ``slightly'' from a fixed generalized flag $W$ in $V$. For the exact definition see Section \ref{section:brief_background}. One then shows that the so obtained ind-variety is isomorphic to $G/P$ for $G=\SL(\infty)$ and some splitting parabolic subgroup $P\subset G$.
	
	An ind-grassmannian is an ind-variety of generalized flags for which the fixed generalized flag consists of a single proper subspace $W\subset V$. For $\dim W=\codim_V W=\infty$ 
	the ind-grassmannian is isomorphic to the Sato grassmannian. This has been pointed out for instance in \cite{GattoSalehyan}.
	
	In the cases of the groups $\SO(\infty)$ and $\Sp(\infty)$ we consider ind-varieties of isotropic generalized flags, as stated in Section \ref{section:more_background}.
	
	Our main result is the explicit determination of the group $\mathrm{Aut}\,X$ for an arbitrary ind-variety of, possibly isotropic, generalized flags. A notable feature is that the answer is very different from the ind-groups $\Pro \GL(\infty)$, $\Pro \Orth(\infty)$, or $\Pro \Sp(\infty)$, and we present it in the language of Mackey groups. Such a group is defined in terms of a nondegenerate pairing of vector spaces $T\times R\rightarrow \Cp$, and is a subgroup of the group of all linear operators $\vfi \colon T\rightarrow T$ for which the dual operator $\vfi^*$ determines a well-defined automorphism $\overline{\vfi}\colon R\rightarrow R$. This definition of Mackey group is inspired by G. Mackey's dissertation \cite{Mackey45}. If $T$ and $R$ are finite dimensional, then the Mackey group is nothing but $\GL(T)\simeq\GL(R)$. The group known as Japanese $GL(\infty)$ is a Mackey group and plays a crusial role in our work. In the Appendix we discuss the structure of this group in detail.
	
	The precise statement of our main result, Theorem \ref{theorem:1}, is presented in Section \ref{section:brief_background}. The consideration of the isotropic case is postponed to Section \ref{sect:proof_theorem_2_isotropic_case}. The proof of Theorem \ref{theorem:1} is divided into two parts: the case of an ind-grassmannian and the case of an arbitrary generalized flag. For the Sato grassmannian (which is the most interesting ind-grassmannian) our result implies that its automorphism group is isomorphic to the projectivization of the connected component of the identity in the group Japanese $GL(\infty)$. In Section \ref{sect_explicit_matrix_form} we give a matrix realization of the group of automorphisms of an arbitrary ind-variety of generalized flags. In the isotropic case such a realization is given in Corollary \ref{corollary:group_structure_isotropic_case}.
	
	
	We would like to point out that $\mathrm{Aut}\, X$ depends essentially on the ind-variety $X$, despite the fact that all $X$ are homogeneous spaces for the same group $\SL(\infty)$ (or, respectively, $\SO(\infty)$, $\Sp(\infty)$). This is in contrast with the finite-dimensional case in which the connected component of the identity in the automorphism group of a variety $\SL(n)/P$ (respectively, $\SO(2n)/P$, $\SO(2n+1)/P$ or $\Sp(2n)/P$) depends only on $n$ and not on the choice of $P$. Further research should be carried out to compare the isomorphism classes of ind-varieties of generalized flags with the isomorphism classes of their automorphism groups.
	
	Our possible application of the results of the present paper would be the study of locally reductive ind-groups $\widetilde{G}$ different from $G=\SL(\infty)$, $\SO(\infty)$, $\Sp(\infty)$ for which $G/P$ is a homogeneous $\widetilde{G}$-space.
	
	In conclusion of this short introduction, we should mention that some particular cases of the automorphism groups of ind-varieties of generalized flags have been considered in \cite{Penkov} and \cite{Tarigradschi}.
	
	\textbf{Acknowledgements}. 
	
	I.P. thanks Francesco Esposito for inspiring discussion on topics related to the present paper. M.I. was supported by the Foundation for the Advancement of Theoretical Physics and Mathematics ``BASIS'' grant no. 18--1--7--2--1, and by RFBR grant no. 20--01--00091--a. I.P. was supported in part by DFG grant PE 980/8--1.
	
	We are grateful to two referees for their sharp comments. One of the referees found several inaccuracies and made helpful suggestions for improving the text. Another referee suggested to mention that in the works \cite{Huckleberry1}, \cite{Huckleberry2} of Alan Huckleberry the classical result of \cite{Onishchik} has been extended to the non-compact case of flag domains.

	\sect{Brief background and statement of the main result}\label{section:brief_background}
	
	The ground field is $\Cp$. If $R$ is a vector space, we set $\GL(R)=\{ \varphi\in \mathrm{Hom}_{\Cp}(R,R)\mid \varphi\text{ is invertible}\}$ and $R^*=\mathrm{Hom}_{\mathbb{C}}(R,\mathbb{C})$. We also use the superscript $\cdot^*$ to denote the dual of a vector bundle, as well as the pullback functor for vector bundles along a morphism of varieties. In what follows we consider infinite matrices, in particular, infinite rows and columns. We call such matrices, rows or columns \emph{finitary} if they have at most finitely many nonzero entries. By $\langle \cdot \rangle_{\Cp}$ we denote the linear span over $\Cp$.
	
	We fix a countable-dimensional vector space $V$. A \emph{chain of subspaces} $W=\{W_{\alpha}\}$ is a set of subspaces $W_{\alpha}\subseteq V$, parameterized by some index set with elements $\alpha$ such that for $\alpha\neq \alpha'$ we have $W_{\alpha}\varsubsetneq W_{\alpha'}$ or $W_{\alpha'}\varsubsetneq W_{\alpha}$. The relation of inclusion induces a total order on the set of indices of a chain. A chain of subspaces  $W=\{W_{\alpha}\}$ is a \emph{generalized flag in} $V$ if every index $\alpha$ has either an immediate predecessor or an immediate successor, and every nonzero vector $v$ of $V$ is contained in some difference $W_{\alpha''}\setminus W_{\alpha'}$, where $\alpha''$ is the immediate successor of $\alpha'$. For a more detailed discussion of generalized flags, and for an introduction to ind-varieties of generalized flags, see, e.g., \cite{DimitrovPenkov}, \cite{IgnatyevPenkov} and \cite{PenkovTikhomirov2}.
	
	We say that a generalized flag $W$ is \emph{compatible with a basis} $\widetilde{E}$ of $V$ (or that $W$ is $\widetilde{E}$-\emph{compatible}) if any space $W_{\alpha}$ of $W$ is spanned by elements of $\widetilde{E}$, i.e., $W_{\alpha}=\langle W_{\alpha}\cap \widetilde{E} \rangle_{\Cp}$ for any $\alpha$. We set $E_{\alpha}=\widetilde{E}\cap W_{\alpha}$. Then $E_\alpha$ is a basis of~$W_{\alpha}$. By $(W_{\alpha})_*$ we denote the span of the system of linear functions $\widetilde{E}^*$ dual to the basis $E_{\alpha}$. We have $(W_{\alpha})_*\subseteq W^*_{\alpha}$. We also let $V_*$ equal the span of the system of linear functions dual to the basis $\widetilde{E}$. The group $\GL(\widetilde{E},V)$ is the subgroup of $\GL(V)$ consisting of all invertible operators $\vfi\colon V\rightarrow V$ each of which acts as the identity on all but finitely elements of $\widetilde{E}$.
	
	In what follows, we fix a basis $\widetilde{E}$ and a generalized flag $W$ in $V$ compatible with $\widetilde{E}$. The set $\Fl(W,\widetilde{E},V)$ is the set of all generalized flags $W'=\{W'_{\alpha}\}$ which are $\widetilde{E}$\emph{-commensurable} with $W$. This latter requirement spells out as the following three conditions for each element $W'$ of $\Fl(W,\widetilde{E},V)$:
	
	\begin{itemize}
		\item the index set which parameterizes the generalized flag $W'$ is the same as the index set of $W$;
		\item there exists a finite-dimensional subspace $Z\subsetneq V$ depending on $W'$, such that for any $\alpha$ we have $W'_{\alpha}+Z= W_{\alpha}+Z$ and $\dim(W_{\alpha}\cap Z)=\dim(W'_{\alpha}\cap Z)$;
		\item the generalized flag $W'$ is compatible with a basis $\widetilde{E}'$ of $V$, depending on $W'$, such that $\widetilde{E}'$ differs from $\widetilde{E}$ by finitely many vectors.
	\end{itemize}
	
	The set $\Fl(W,\widetilde{E},V)$ has a natural structure of ind-variety. This is explained in detail in \cite{DimitrovPenkov} (and in \cite{IgnatyevPenkov} and \cite{PenkovTikhomirov2}). Briefly, the ind-variety structure on $\Fl(W,\widetilde{E},V)$ arises as follows. Enumerate the basis $\widetilde{E}$ by the set $\Zp_{>0}$ and put $V_{i}:=\langle e_1,\ldots,e_i \rangle_{\Cp}$ for $i\in\Zp_{>0}$. Each intersection $W\cap V_i$ is a flag in $V_i$ of certain type $\underline{d}_i=(d_i^1,\ldots,d_i^{k_i})$, and this ordering of the basis $\widetilde{E}$ induces embeddings \begin{equation}\label{formula:flag_embeddings}
		\Fl(\underline{d}_i,V_i)\hooklongrightarrow \Fl(\underline{d}_{i+1},V_{i+1}),
	\end{equation}
	called \emph{strict standard extensions},
	such that $\Fl(W,\widetilde{E},V)=\ilm \Fl(\underline{d}_i,V_i)$. The embeddings (\ref{formula:flag_embeddings}) endow $\Fl(W,\widetilde{E},V)$ with an ind-variety structure. In Section \ref{section:more_background} below we recall the definition of a strict standard extension.
	
	Next we recall that if $T$ and $R$ are two (in general, infinite-dimensional) vector spaces endowed with a non-degenerate pairing $\mathbf{p}\colon T\times R \rightarrow \Cp$, then the Mackey group $ G(T,R)$ is defined as \begin{equation}\label{formula:first_definition_mackey_group}
		G(T,R)=\{\varphi\in\GL(T)\mid \varphi^*(R)= R\},
	\end{equation}
	see \cite{Mackey45}. Here $\varphi^*\colon T^*\rightarrow T^*$ is the operator dual (adjoint) to the operator $\varphi\colon T\rightarrow T$, and $R$ is considered as a subspace of $T^*$ via the embedding $R\hookrightarrow T^*$ induced by the pairing $\mathbf{p}$. Equivalently, $ G(T,R)$ can be defined as the group
	\begin{equation}\label{formula:second_definition_mackey_group}
		\{\psi\in\GL(R)\mid \psi^*(T)= T\}
	\end{equation}
	where $T$ is considered as a subspace of $R^*$ via $\mathbf{p}$. The correspondence $$\varphi\mapsto \left(\varphi^*|_R\right)^{-1}$$ is a canonical isomorphism between the groups (\ref{formula:first_definition_mackey_group}) and (\ref{formula:second_definition_mackey_group}). In what follows, when writing $\vfi \in  G(T,R)$ we will assume that $\vfi\in \GL(T)$, and will denote the operator $\vfi^*|_R$ by $\overline{\vfi}$. Note that, given a subspace $A$ of~$T$, one has $\vfi(A)=\overline{\vfi}^{-1}(A^{\perp})$, where $A^{\perp}$ is the annihilator of $A$ in $R$.
	
	Consider again the ind-variety of generalized flags $\Fl(W,\widetilde{E},V)$. Define the spaces $V^W_{\widetilde{E}}$ and $V^W_{*\widetilde{E}}$ as \begin{equation*}
		\begin{split}	
			V^W_{\widetilde{E}}&:=\bigcap\limits_{\alpha} \left(\left(\left(W_{\alpha}\right)_*\right)^*\oplus U_{\alpha}\right),\\
			V^W_{*\widetilde{E}}&:=\bigcap\limits_{\alpha} \left(\left(W_{\alpha}\right)_*\oplus U^*_{\alpha}\right),
		\end{split}
	\end{equation*}
	where $\alpha$ runs over the indices parameterizing the generalized flag $W$, and the spaces $U_{\alpha}$ are direct complements of the spaces $W_{\alpha}$, i.e., $V=W_{\alpha}\oplus U_{\alpha}$, with the assumption that $U_{\alpha}\cap \widetilde{E}$ is a basis of $U_{\alpha}$. 
	Note that the spaces $V^W_{\widetilde E}$, $V^W_{*\widetilde E}$ are necessarily infinite-dimensional and there is a canonical non-degenerate pairing \begin{equation*}
		V^W_{\widetilde{E}}\times V^W_{*\widetilde{E}} \rightarrow \Cp,
	\end{equation*}  
	therefore the group $ G(V^W_{\widetilde{E}},V^W_{*\widetilde{E}})$ is well defined.
	
	An essential observation is that the spaces $V^W_{\widetilde{E}}$ and $V^W_{*\widetilde{E}}$ depend only on the ind-variety $\Fl(W,\widetilde{E},V)$ and not on the specific point $W\in \Fl(W,\widetilde{E},V)$. This follows from the fact that, for each $\alpha$, the spaces $((W_{\alpha})_*)^*\oplus U_{\alpha}$ and $(W_{\alpha})_*\oplus U^*_{\alpha}$ do not change when $W_{\alpha}$ is replaced by a subspace $W_{\alpha}'\subset W$ which is $\widetilde{E}$-commensurable with $W_{\alpha}$, and $U_{\alpha}$ is replaced by a direct complement $U_{\alpha}'$ of $W_{\alpha}'$ containing all but finitely many vectors from $\widetilde{E}\cap U_{\alpha}$.

	Moreover, $\GL(\widetilde{E},V)$ is a subgroup of $ G(V^W_{\widetilde{E}},V^W_{*\widetilde{E}})$. To see this, consider a linear operator $\varkappa:V\to V$, $\varkappa\in GL(\widetilde E,V)$, and fix $\alpha$. There exist subspaces $W_{\alpha}'\subset W_{\alpha}$, $U_{\alpha}'\subset U_{\alpha}$, such that $\restr{\varkappa}{W_{\alpha}'}=\id_{W_{\alpha}'}$, $\restr{\varkappa}{U_{\alpha}'}=\id_{U_{\alpha}'}$, and $V=W_{\alpha}'\oplus K\oplus U_{\alpha}'$ for some $\varkappa$-invariant finite-dimensional subspace $K\subset V$. Then $(W_{\alpha}')_*$, $(U_{\alpha}')^*$ and $K^*$ are $\varkappa^*$-invariant subspaces of $V^*$, and $(W_{\alpha})_*\oplus U_{\alpha}^*=(W_{\alpha}')_*\oplus K^*\oplus(U_{\alpha}')^*$ is a $\varkappa^*$-invariant subspace of $V^*$. This shows that $V^W_{*\widetilde E}=\bigcap_{\alpha}((W_{\alpha})_*\oplus U_{\alpha}^*)$ is a $\varkappa^*$-invariant subspace of $V^*$. The same argument applied to the $\varkappa^*$-invariant subspace $V_*$ of $V^*$ implies that $V^W_{\widetilde E}$ is a $(\restr{\varkappa^*}{V_*})^*$-invariant subspace of $(V_*)^*$. This allows to consider~$\varkappa$ as an element of $ G(V^W_{\widetilde{E}},V^W_{*\widetilde{E}})$.

	Next, if $W'$ is any chain of subspaces in $V$ and $\varphi\colon V^W_{\widetilde{E}}\rightarrow V^W_{\widetilde{E}}$ is any linear operator from the group $ G(V^W_{\widetilde{E}},V^W_{*\widetilde{E}})$, then \begin{equation}\label{formula:action_on_flag}
		\overline{\varphi}^{-1}\left({W'}^\perp\right)^{\perp}\cap V=\vfi((W'^{\perp})^{\perp})\cap V
	\end{equation}
	is a chain of subspaces in $V$. Here $W'^{\perp}$ is the chain in $V^W_{*\widetilde{E}}$ consisting of the annihilators in $V^W_{*\widetilde{E}}\subseteq V^*$ of the spaces $W_{\alpha}'$, and similarly $\overline{\varphi}^{-1}\left(W'^\perp\right)^{\perp}$, $W'^{\perp\perp}$ are chains in $V^W_{\widetilde{E}}\subseteq V^{**}$. 
	
	In what follows, we use the notation ${W'}^{\perp}$ for chains perpendicular to $W'$ also in appropriate subspaces of $V^*$ different from $V^W_{*\widetilde{E}}$, and indicate the respective subspace as necessary. A similar convention applies to the notation ${W'}^{\perp \perp}$. Moreover, we call the generalized flag $W$ \emph{symmetric} if the chain $W^{\perp}\subset V_*$ is the image of $W$ under a linear isomorphism $V_*\simeq V$ sending $\widetilde{E}^*$ to $\widetilde{E}$.
	

	\theopwop{\label{theorem:1}
		\begin{itemize}
			\item[\textup{a)}] If $W$ is not symmetric\textup, then the group $\mathrm{Aut}\,\Fl(W,\widetilde{E},V)$ is isomorphic to $\Pro(\GL(\widetilde{E},V) \cdot \St_W )$\textup, where $\St_W$ is the stabilizer of the generalized flag $W$ in the group $ G(V^W_{\widetilde{E}},V^W_{*\widetilde{E}})$ under the action~\textup{(\ref{formula:action_on_flag})}. Here the product $\cdot$ is taken inside $ G(V^W_{\widetilde{E}},V^W_{*\widetilde{E}})$\textup, and $\Pro\cdot$ indicates passage to the quotient modulo scalar operators.
			\item[\textup{b)}] If $W$ is symmetric\textup, then the group $\mathrm{Aut}\,\Fl(W,\widetilde{E},V)$ is isomorphic to $\Pro(\GL(\widetilde{E},V) \cdot \St_W)\rtimes \mathbb{Z}_2$.
	\end{itemize}}
	
	In Section \ref{sect_explicit_matrix_form} we present an explicit matrix realization of the group $\GL(\widetilde{E},V) \cdot \St_W$. Let's also point out that, since Theorem~\ref{theorem:1} implies that $\GL(\widetilde{E},V) \cdot \St_W$ is a group, we have $\GL(\widetilde{E},V) \cdot \St_W=\St_W\cdot \GL(\widetilde{E},V)$.
	
	\nota{In the case of a finite-dimensional flag variety $X$, every automorphism of $X$ belonging to the connected component of unity in the automorphism group has a fixed point on $X$. This no longer holds in the generality of Theorem \ref{theorem:1}. Indeed, if $X$ is the projective ind-space $\Fl(W,E,V)$ for $\dim W=1$, then $\mathrm{Aut}X=\PGL(V)$ and it is well known that not every invertible linear automorphism of $V$ has an eigenvector.}

	\newpage\sect{Examples}
	
	Before we embark on proving Theorem \ref{theorem:1}, we present five examples in which we compute the respective group $\GL(\widetilde{E},V) \cdot \St_W$ from Theorem \ref{theorem:1}. 
	In all five cases our claims follow from Theorem \ref{theorem:description_general_case} below, which provides a matrix form of the group $\GL(\widetilde{E},V) \cdot \St_W$ in the general case.
	
	\sst{The case of an ind-grassmannian} Let's consider the case where the generalized flag $W$ has the form $0\subsetneq W\subsetneq V$, where $W$ is a single proper subspace of $V$ (we slightly abuse notation by using the same letter $W$ for a flag and a subspace). There are three cases: $\dim W<\infty$, or $\dim W=\codim_V W=\infty$, or $\codim_V W<\infty$. If $\dim W<\infty$, then $\Gra(W,\widetilde{E},V)$ does not depend on the basis $\widetilde{E}$, and the points of $\Gra(W,\widetilde{E},V)$ are all subspaces of $V$ of the same dimension as $W$. In this case we may write $\Gra(W,\widetilde{E},V)=\Gra(\dim W, V)$. If $\codim_V W<\infty$ then $\Gra(W,\widetilde{E},V)$ depends as a set on the choice of the basis $\widetilde{E}$, but up to isomorphism $\Gra(W,\widetilde{E},V)$ depends only on $\codim_V W$. Moreover, as an ind-variety $\Gra(W,\widetilde{E},V)$ is isomorphic to $\Gra(\codim_V W,V)$. If $\dim W=\codim_V W=\infty$, the ind-variety $\Gra(W,\widetilde{E}, V)$ does not depend up to isomorphism on the choice of both $W$ and $\widetilde{E}$. It is known, see for instance \cite{GattoSalehyan}, that in this case $\Gra(W,\widetilde{E}, V)$ is isomorphic to the Sato grassmannian introduced in \cite{Sato}.
	
	If $\dim W<\infty$, then $V^W_{\widetilde{E}}=V$, $V^W_{*\widetilde{E}}=V^*$ and $ G(V,V^*)=\GL(V)$. Since $W$ is not symmetric, Theorem \ref{theorem:1} asserts that $\mathrm{Aut}\,\Gra(W,\widetilde{E},V)\cong \Pro(\GL(\widetilde{E},V) \cdot \St_W)$, and we note that here  \begin{equation}\label{formula:automorphism_group_grassmannian_special_case}
		\Pro(\GL(\widetilde{E},V) \cdot \St_W)\cong\Pro \GL(V).
	\end{equation}
	Indeed, the action of $ G(V,V^*)=\GL(V)$ on $\Gra(W,\widetilde{E},V)$ via the formula (\ref{formula:action_on_flag}) is easily checked to coincide with the obvious action of $\GL(V)$ on subspaces of $V$, and the isomorphism (\ref{formula:automorphism_group_grassmannian_special_case}) is a consequence of the transitive action of the group $\GL(\widetilde{E},V)$ on finite-dimensional subspaces of fixed dimension in~$V$.
	
	If $\codim_V W<\infty$, then $V^W_{\widetilde{E}}=(V_*)^*$, $V^W_{*\widetilde{E}}=V_*$ and  $\Pro(\GL(\widetilde{E},V) \cdot \St_W)\cong\Pro \GL(V_*)$, i.e.,
	\begin{equation*}
		\mathrm{Aut}\,\Gra(W,\widetilde{E},V)\cong\Pro \GL(V_*).
	\end{equation*}
	
	
	In the case where $\dim W=\codim_V W=\infty$, we prove in the Appendix that the group $G(V^W_{\widetilde{E}},V^W_{*\widetilde{E}})$ can be represented as invertible $(\Zp\setminus\{0\})\times(\Zp\setminus\{0\})$-matrices which together with their inverses satisfy the condition: in the block structure 
	\begin{equation}\label{formula:block_matrix_structure}
		\left(\begin{array}{c|c}
			A & B \\
			\hline
			C & D  \end{array}\right)
	\end{equation}
	induced by the equality $\Zp\setminus\{0\}=\Zp_{<0}\sqcup\Zp_{>0}$, the matrix $A$ has finitary rows (no restriction on the columns), the matrix $D$ has finitary columns (no restriction on the rows) and the matrix $C$ is finitary. 
	The group $\GL(\widetilde{E},V) \cdot \St_W$ consists of matrices $M$ such that $M$ and $M^{-1}$ have the form (\ref{formula:block_matrix_structure}) and satisfy the additional condition $\rk C=\rk C'$ where \begin{equation*}
		M^{-1}=\left(\begin{array}{c|c}
			A' & B' \\
			\hline
			C' & D' \end{array}\right).
	\end{equation*}  Moreover, in this case $W$ is symmetric.  
	
	
	\sst{The case of $\Fl(W,\widetilde{E},V)$, where $W=\{W_n\}$, $\dim W_n=n$ for $n\in \Zp_{>0}$, $\bigcup_n W_n=V$}\label{example:generalized_flag_full} In this case $V^W_{\widetilde{E}}=V$, $V^W_{*\widetilde{E}}=V^*$, and the group $\GL(\widetilde{E},V) \cdot \St_W$ can be identified with all invertible  $\Zp_{>0}\times \Zp_{>0}$-matrices with finitely many nonzero entries below the main diagonal, cf. \cite{Penkov}. 
	
	\sst{The case of $\Fl(W,\widetilde{E},V)$, where $\widetilde{E}=\{e_i\}_{i\in \Zp}$ and $W=\{W_n=\langle e_i,~i\leq n \rangle_{\Cp},~n\in \Zp\}$}\label{example:generalized_flag_basis_E_tilde} Here $$V^W_{\widetilde{E}}=\left(\left(W_0\right)_*\right)^*\oplus  \langle e_i,~i>0\rangle_{\Cp},~V^W_{*\widetilde{E}}=(W_0)_*\oplus \left( \langle e_i,~i>0\rangle_{\Cp}\right)^*.$$ In coordinate form, the vectors from~$V^W_{\widetilde{E}}$ are columns $(a_j)_{j\in \Zp}$ with $a_j=0$ for $j\gg 0$, and 
	$\GL(\widetilde{E},V) \cdot \St_W$ consists of   all invertible $\Zp  \times \Zp $-matrices $M$ which, together with their inverses, have finitely many nonzero entries below the main diagonal and satisfy the condition $\rk C=\rk C'$, where $C$ and $C'$ are respectively the strictly lower-triangular parts of $M$ and $M^{-1}$.
	
	\sst{The case of $\Fl(W,\widetilde{E},V)$, where $W=\{0\subset W_1 \subset W_{-1} \subset V \}$, $\dim W_1=1$, $\codim_V W_{-1}=1$}\label{example:generalized_flag_W_1_W_-1} Here $\widetilde{E}$ can be ordered by any countable ordered set $I$ with a minimal and a maximal element. We have $V^W_{\widetilde{E }}=V$, $V^W_{*\widetilde{E}}=V_*$, and $\GL(\widetilde{E},V) \cdot \St_W$ consists of all invertible $I\times I$-matrices 
	which, together with their inverses, satisfy the condition that each row and column is finitary. 
	
	\sst{The case of $\Fl(W,\widetilde{E},V)$ for $\widetilde{E}=\{e_i\}_{i\in \Zp_{>0} \sqcup \Zp_{<0}}$ and $W=\{ W_n \}$, where \linebreak $W_n=\langle e_1,\ldots,e_n \rangle_{\Cp}$ if $n\in \Zp_{>0}$ and $W_n=\langle\ldots,e_{n-2},e_{n-1},e_1,e_2,\ldots\rangle_{\Cp}$ if $n\in \Zp_{<0}$}\label{example:generalized_flag_basis_E_tilde_pm} Then $V^W_{\widetilde{E}}=V$, $V^W_{*\widetilde{E}}=V_*$, 
	$I=\Zp_{>0} \sqcup \Zp_{<0}$ with $k<l$ if $k\in \mathbb{Z}_{>0}$, $l\in \mathbb{Z}_{<0}$, and $\GL(\widetilde{E},V) \cdot \St_W$ consists of all invertible $I\times I$-matrices which, together with their inverses, satisfy the condition: each row and column 
	is finitary, and there are at most finitely many nonzero entries below the main diagonal.
	\vspace{0.4cm}

	Note that in the cases \ref{example:generalized_flag_basis_E_tilde}, \ref{example:generalized_flag_W_1_W_-1}, and \ref{example:generalized_flag_basis_E_tilde_pm} the generalized flag $W$ is symmetric, while in the case \ref{example:generalized_flag_full} $W$ is not symmetric.

	\sect{More background}\label{section:more_background}

	We need to recall some facts about linear embeddings of finite-dimensional grassmannians and flag varieties. If $T$ is a finite-dimensional space and $\underline{d}=\{d^1,\ldots,d^i\}$ is a vector of positive integers satisfying $d^k<d^l<\dim T$ for $k<l$, then $\Fl(\underline{d},T)$ denotes the variety of all flags of subspaces $T_1\subsetneq \ldots \subsetneq T_i \subsetneq T$, where $\dim T_j=d^j$. If $\underline{d}$ consists of one integer $d$, we write simply $\Gra(d,T)$. If $T$ is endowed with a non-degenerate symmetric or antisymmetric (symplectic) form, we write respectively $\FlO(\underline{d},T)$ and $\FlS(\underline{d},T)$ for the varieties of isotropic flags in $T$ with respect to the fixed form. We also write $\GraO(d,T)$ and $\GraS(d,T)$. 
	An isotropic flag has always length less or equal $\frac{\dim T}{2}$ but, for convenience, in this paper by an \emph{isotropic flag} we will mean a flag of the form $$W_1\subset W_2 \subset \ldots \subset W_k \subset W^{\perp}_k \subset \ldots \subset W^{\perp}_1,$$ where the spaces $W_1,~\ldots,~W_k$ are isotropic and the spaces $W^{\perp}_k,~\ldots,~W^{\perp}_1$ are coisotropic. All flag varieties $\FlO(\underline{d},T)$ and $\FlS(\underline{d},T)$ are connected, except 
	$\GraO(d,T)$ for $\dim T=2d$. In what follows, by $\FlO(\underline{d},T)$ or $\GraO(\underline{d},T)$ we always denote a connected component.
	
	The Picard group of any grassmannian or ind-grassmannian $Z$ is isomorphic to $\Zp$ except in the case of $\GraO\left(\frac{\dim T}{2}-1,~T\right)$ for $2<\dim T\in 2\Zp_{>0}$, and $\mathcal{O}_Z(1)$ always denotes the ample generator of $\mathrm{Pic}\,Z$. In the case of $\GraO\left(\frac{\dim T}{2}-1,~T\right)$ for $2<\dim T\in 2\Zp_{>0}$ we have $\mathrm{Pic}\, \GraO\left(\frac{\dim T}{2}-1,~T\right)\simeq \Zp \times \Zp$.  
	
	The automorphism groups of the flag varieties $\Fl(\underline{d},~T)$, $\FlO(\underline{d},~T)$, $\FlS(\underline{d},~T)$ have been known for long time. The fact that the automorphism group of the projective space $\mathbb{P}^n$ is $\PGL (n+1)$ goes back to the nineteenth century. Wei-Liang Chow \cite{Chow} extended this result to grassmannians in 1949. 
	For a general flag variety $X=\Fl(\underline{d},~T)$, $\FlO(\underline{d},~T)$, $\FlS(\underline{d},~T)$ 
	the connected component of the identity in the automorphism group $\mathrm{Aut}\,X$ is the respective group $\PGL (T)$, $\SO (T)$, or $\Sp (T)$, except in several cases listed by A.L. Onishchik in \cite{Onishchik}. These special cases are 
	$\GraS(1,T)\simeq\mathbb{P}(T)$, the five dimensional quadric $\GraO(1,T)$ for $\dim T=7$, and $\GraO \left(\frac{\dim T -1}{2},~T\right)$ for $\dim T\in 2\Zp_{>0}+1$. 
	
	In all cases, see for instance \cite[Section 3.3]{Akhiezer}, the full automorphism group $G$ is always a semidirect product of its connected component of unity $G^0$ and a finite group 
	of automorphisms of the Dynkin diagram  of the Lie algebra $\mathfrak{g}=\mathrm{Lie}\, G^0$ which keep fixed 
	the simple roots of the Lie algebra of the isotropy subgroup of a point on the respective flag variety. In the present paper we only consider classical groups of large enough rank, hence we can summarize the relevant part of this result as follows:
	\begin{itemize}
		\item for $\dim T\geq 3$, 
	\end{itemize}
	$$\mathrm{Aut}\, \Fl (\underline{d},T)\simeq \begin{cases} \Pro \GL(T) \rtimes \Zp_2 \text{ if } \underline{d}=(d^0=0,~d^1,~d^2,~\ldots,~d^s,~d_n^{s+1}=\dim T)\text{ satisfies  }  
		\\ \text{ the condition } d^i-d^{i-1}=d^{s+2-i}-d^{s+1-i}\text{ for all } 1\leq s\leq n+1
		\\
		\PGL(T) \text{ in all other cases. } \end{cases}$$ 
	\begin{itemize}
		\item for $\dim T\geq 8$, $\mathrm{Aut}\, \FlO (\underline{d},T)\simeq \Orth (T)$, except for $\GraO(\frac{\dim T}{2},T)$ where $\mathrm{Aut}\, \GraO(\frac{\dim T}{2},T)\simeq \SO (T)$, and for $\GraO(\frac{\dim T-1}{2},T)$ where $\mathrm{Aut}\, \GraO(\frac{\dim T-1}{2},T)\simeq \SO (T')$ for $\dim T'=\dim T +1$.
		\item for $\dim T\geq 4$, $\mathrm{Aut}\, \FlS(\underline{d},T)\simeq \Sp(T)$, except for $\GraS(1,T)$ where  $\mathrm{Aut}\,\GraS(1,T)\simeq \Pro \GL(T)$.
	\end{itemize}
	Based on the above exception concerning $\GraO\left(\frac{\dim T -1}{2},~T\right)$ for $\dim T\in 2\Zp_{>0}+1$, in what follows we will automatically assume that this case is excluded from consideration. This leads to no loss of generality as $\GraO\left(\frac{\dim T -1}{2},~T\right)$ is isomorphic to $\GraO\left(\frac{\dim T'}{2},~T'\right)$ where $T'$ is an orthogonal space of dimension $\dim T+1$.
	
	A nice class of embeddings of flag varieties 
	$\Fl(\underline{d}_1,T_1)\hookrightarrow \Fl(\underline{d}_2,T_2)$ for $\dim T_1<\dim T_2$ is the class of standard extensions. Embeddings of ind-grassmannians are discussed in detail in \cite{PenkovTikhomirov1}, and of arbitrary flag varieties in \cite{PenkovTikhomirov2}. 
	Here we just recall a definition and a basic fact needed to understand our arguments in Section \ref{sect:proof_theorem_1_for_ind-grassmannians}, \ref{sect:proof_theorem_1_in_the_general_case}, \ref{sect:proof_theorem_2_isotropic_case}.
	
	\defi{\label{descrn of st eTt}\begin{itemize}
			\item[a)] An embedding \begin{equation*}
				\eta \colon \Fl(d^1_1,\ldots,d^k_1,T)\hookrightarrow 
				\Fl(d^1_2,\ldots,d^l_2,T'),
			\end{equation*} respectively, \begin{equation*}\eta\colon \FlO(d^1_1,\ldots,d^k_1,T)\hookrightarrow \FlO(d^1_2,\ldots,d^l_2,T'),\end{equation*} respectively, \begin{equation*}\eta\colon \FlS(d^1_1,\ldots,d^k_1,T)
				\hookrightarrow \FlS(d^1_2,\ldots,d^l_2,T')\end{equation*} 
			is a \emph{strict standard extension} if there exists a surjection
			\begin{equation*}p\colon\{0,1,\ldots,l,l+1\}\to
				\{0,1,\ldots,k,k+1\}\end{equation*} 
			satisfying $p(i)\leq p(j)$ for $i<j$, together with an isomorphism
			\begin{equation}\label{direct sum with W}
				V'=V\oplus \widehat{W},
			\end{equation}
			satisfying $\widehat{W}=V^{\bot}$ in the orthogonal and symplectic 
			case, 
			and subspaces $W_i\in\wh W$ for $1\leq i\leq l$ with $W_i\subset W_j$ for $i<j$ such that $\eta$ has the form
			\begin{equation}\label{phi(...)1}
				\eta\big(\{0\}=V_{p(0)}\subset V_{p(1)}\subset\ldots\subset 
				V_{p(l)}\subset V\big)=
				\big(\{0\}=V_{p(0)}\subset V_{p(1)}\oplus W_1\subset\ldots\subset 
				V_{p(l)}\oplus W_{l}\subset V'\big).
			\end{equation}
			Here the spaces $W_1,\ldots,W_l$ are not required to be pairwise distinct, while the spaces in the right-hand side of (\ref{phi(...)1}) are pairwise distinct by definition. In the orthogonal or symplectic case we require that for each $i$ there is $j$ such that $W_j=W_i^{\perp}$.
			\item[b)] An embedding $$\eta\colon \Fl(d^1_1,~\ldots,~d^k_1,~T)\hooklongrightarrow \Fl(d^1_2,~\ldots,~d^l_2,~T')$$
			is a \emph{standard extension} if after composing with one of the duality isomorphisms
			$$\Fl(d^1_1,~\ldots,~d^k_1,~T)\simeq \Fl(\dim T -d^k_1,~\ldots,~\dim T -d^1_1,~T^* ),$$
			$$\Fl(d^1_2,~\ldots,~d^l_2,~T')\simeq \Fl(\dim T' -d^l_2,~\ldots,~\dim T' -d^1_2,~{T'}^* )$$
			$\eta$ becomes a strict standard extension. For varieties of isotropic flags, standard extension and strict standard extension are synonyms.
	\end{itemize}}
	
	The following theorem follows directly from Corollary 4.4 
	in \cite{PenkovTikhomirov2}.
	
	\theopwop{\label{theorem:P}
		Let $X_1\hookrightarrow X_2$ and $Y_1\hookrightarrow Y_2$ be embeddings of flag varieties or of isotropic flag varieties, such that $Y_1$ is the image of $X_1$ under some isomorphism $\vfi\colon X_2 \rightarrow Y_2$. Then the embedding $Y_1\hookrightarrow Y_2$ is a standard extension whenever the embedding $X_1\hookrightarrow X_2$ is a standard extension.}

	We need to recall also the notion of an ind-variety of isotropic generalized flags. There are several cases. If a symmetric non-degenerate form $(\cdot,\cdot)$ on $V$ is given, then there are two types of relevant bases $\widetilde{E}$ we consider: $$\{ e_i \}_{i\in \mathbb{Z}}\text{ with } (e_i,~e_{-i})=1\text{ for  all }i,~(e_i,~e_k)=0\text{ for }k\neq -i,$$ 
	or \begin{equation}\label{formula:second_basis}
		\{e_i\}_{i\in \mathbb{Z}\setminus\{0\}}\text{ with }(e_i,~e_{-i})=1 \text{ for all }i\in \mathbb{Z}_{>0},~ (e_i,~e_k)=0\text{ for }k\neq -i.
	\end{equation}
	In the case of a symplectic non-degenerate form on $V$, we consider bases satisfying (\ref{formula:second_basis}) (here $(e_{-i},e_i)=-1$ for $i\in \mathbb{Z}_{>0}$). 
	We refer to bases as above as \emph{isotropic bases} of $V$. An \emph{isotropic generalized flag} in $V$ is by definition a generalized flag $W=\{ W_{\alpha} \}$ in $V$ such that each space $W_{\alpha}$ is either isotropic or coisotropic and  
	$W_{\alpha}$ belongs to $W$ if and only if $W^{\perp}_{\alpha}=\{ w'\in V\mid (w',w)=0\text{ for all }w\in W_{\alpha} \}$ belongs to $W$. If $W$ is an $E$-compatible isotropic generalized flag in $V$, then by  definition, $\FlO (W,\widetilde{E},V)$ in the case of a symmetric form, or $\FlS(W,\widetilde{E},V)$ in the case of a symplectic form, consists of all isotropic generalized flags in $V$ which are $\widetilde{E}$-commensurable with $W$. 
	In all cases, $\FlO(W,\widetilde{E},V)$ or, respectively, $\FlS(W,\widetilde{E},V)$ is a direct limit of finite-dimensional varieties of isotropic flags under standard extensions.
	
	Finally, we should point out that in order to follow the proof of our main results in Sections \ref{sect:proof_theorem_1_for_ind-grassmannians}--\ref{sect:proof_theorem_2_isotropic_case} the readers should first familiarize themselves with the results of the Appendix.

	\sect{Proof of Theorem \ref{theorem:1} for ind-grassmannians}\label{sect:proof_theorem_1_for_ind-grassmannians}
	
	We start by proving Theorem \ref{theorem:1} under the assumption that 
	the generalized flag $W$ has exactly one proper subspace, which we also denote by $W$. In what follows, we write most of the time $\Gra(W,E,V)$ instead of $\Gra(W,\widetilde{E},V)$, where $E=\widetilde{E}\cap W$. The set $E$ is a basis of $W$, and the ind-variety $\Gra(W,\widetilde{E},V)$ depends only on $E$ and not on the entire basis $\widetilde{E}$. We feel that this notation makes the argument more transparent. Also, the space $W$ is fixed and we write $V_E$ and $V_{*E}$ instead of $V^W_{\widetilde{E}}$ and $V^W_{*\widetilde{E}}$, respectively.
	
	Recall that $\Gra(W,E,V)$ is defined as the direct limit of strict standard extensions \begin{equation*}
		\Gra(d_n,V_n) \hooklongrightarrow \Gra(d_{n+1},V_{n+1})
	\end{equation*}
	for some $d_n$ and some subspaces $V_n\subset V$, $\dim V_n\gee n$, $\ilm V_n=V$. If $d_n$ stabilizes at $k\in\mathbb{Z}_{>0}$ for large $n$, then $\ilm \Gra(d_n,V_n)=\Gra(k,V)$ is the ind-grassmannian of all $k$-dimensional subspaces in $V$. If $\dim{V_n}-d_n$ stabilizes at $k>0$, then $\Gra(W,E,V)$ is isomorphic to $\Gra(k,V)$ as an ind-variety via the map
	\begin{equation*}
		\delta \colon \Gra(W,E,V) \rightarrow \Gra(k,V_*),
	\end{equation*}
	\begin{equation*}
		W'\longmapsto W'^{\perp}\subset V_*,
	\end{equation*}
	where $W'$ denotes a variable point of $\Gra(W,E,V)$ and ${W'}^{\perp}:=\{ \alpha \in V_*\mid 
	\alpha(w')=0~\forall~w' \in W' \}$.
	
	As we already mentioned, the automorphism groups of (finite-dimensional) grassmannians have been described in the classical paper \cite{Chow}. This description implies that if
	\begin{equation}\label{formula:isomorphism_of_grassmannians}
		\alpha \colon \Gra(d_n,T_n)\stackrel{\sim}{\rightarrow}
		\Gra(d_n,T'_n)
	\end{equation}
	is any isomorphism of grassmannians, where $\dim T_n=\dim T'_n$, then the pullback $\alpha^*S'_n$ of the tautological bundle $S'_n$ on $\Gra(d_n,T'_n)$ is isomorphic to the tautological bundle $S_n$ on $\Gra(d_n,T_n)$, or to the bundle $(\widetilde{T}_n/S_n)^*$ in case $2d_n=n$, where $\widetilde{T}_n$ is the trivial bundle on $\Gra(d_n,T_n)$ with fiber $T_n$. Moreover, if $\alpha^*S'_n\simeq S_n$, then the isomorphism (\ref{formula:isomorphism_of_grassmannians}) is determined by the linear operator $\eta \colon(T'_n)^*\rightarrow T^*_n$ which it induces via pullback: we have $\alpha(T_{d_n})=\eta^*(T_{d_n})$ where $T_{d_n}\in \Gra(d_n,T_n)$ and the operator $\eta^*$ is dual to $\eta$. Recall also that any global endomorphism of the bundle $S_n$ or $(\widetilde{T}/S_n)^*$ is scalar. 
	
	Set $X_n=\Gra(d_n,V_n)$. Then $\Gra(W,E,V)=\ilm X_n$. In the rest of the argument we assume in addition that $\dim W=\codim_V W=\infty$ and that $d_n=n$, $\dim V_n=2n$. This is the case of the Sato grassmannian. The remaining cases where $\dim W<\infty$ or $\codim_V W<\infty$ have been considered in \cite{Penkov}, and it has been proved there that $\mathrm{Aut}\,\Gra(W,E,V)\cong\Pro \GL(V)$ for $\dim W<\infty$, and $\mathrm{Aut}\,\Gra(W,E,V)\cong\Pro \GL(V_*)$ for $\codim W<\infty$. This is in agreement with Theorem \ref{theorem:1}, as in these two cases $W$ is not symmetric and there is an isomorphism \begin{equation}\label{formula:group_special_cases}
		\GL(\widetilde{E},V)\cdot \St_W \cong\begin{cases} \GL(V) \text{ for } \dim W<\infty \\ \GL(V_*) \text{ for } \codim_V W<\infty.\end{cases}
	\end{equation} 
	We leave the proof of (\ref{formula:group_special_cases}) as an exercise to the reader. 
	
	Our first step will be to prove that the group of automorphisms of the ind-variety $\Gra(W,E,V)$ is a subgroup of $\Pro(\GL(\widetilde{E},V)\cdot \St_W)\rtimes \mathbb{Z}_2$. Let $\widetilde{\varphi}\colon \Gra(W,E,V)\rightarrow \Gra(W,E,V)$ be an arbitrary automorphism. Denote by $Y_n$ the image of $X_n$, that is, $\widetilde{\varphi}(X_n)=Y_n$. Clearly, $Y_n$ is a grassmannian isomorphic to $X_n$. Moreover, the embeddings $Y_n\hookrightarrow Y_{n+1}$ can be assumed to be strict standard extensions by Theorem \ref{theorem:P}. 
	
	Next, we have two possibilities: for some $n$, the isomorphism \begin{equation*}
		\widetilde{\varphi}_n:=\widetilde{\vfi}|_{X_n}\colon X_n\longrightarrow Y_n
	\end{equation*}
	has the property \begin{equation}\label{formula:first_possibillity_for_varphi_tilde}
		\widetilde{\varphi}^*_nS_{Y_n}\cong \left(\widetilde{V}_n/S_n\right)^*,
	\end{equation}
	where $S_n$ is the tautological bundle on $X_n$ and $S_{Y_n}$ is the tautological bundle on $Y_n$, or the property \begin{equation}\label{formula:second_possibillity_for_varphi_tilde}
		\widetilde{\varphi}^*_nS_{Y_n}\cong S_n
	\end{equation}
	for all $n$. Since our chains of embeddings $X_n\hookrightarrow X_{n+1}$ and $Y_n\hookrightarrow Y_{n+1}$ are strict standard extensions, if the isomorphism (\ref{formula:first_possibillity_for_varphi_tilde}) holds for some $n$, it must hold for all $n$. However, if this happens we can compose our automorphism $\widetilde{\varphi}$ with the following automorphism which represents an element of $\mathbb{Z}_2$ in the semidirect product from the statement of the theorem:\begin{equation*}
		\delta:\Gra(W,\widetilde{E},V)\longmapsto \Gra(W^{\perp},\widetilde{E}^*,V_*)\longrightarrow \Gra(W,\widetilde{E},V),
	\end{equation*}
	where the left arrow sends $W$ to $W^{\perp}\subset V_*$ and the right arrow is induced by an appropriate linear isomorphism \begin{equation*}
		V_*\longrightarrow V
	\end{equation*} 
	which maps $W^{\perp}$ to $W$. Option (\ref{formula:second_possibillity_for_varphi_tilde}) certainly holds for the composition $\delta \circ \widetilde{\varphi}$, so without loss of generality we can assume in the rest of the argument that (\ref{formula:second_possibillity_for_varphi_tilde}) holds.
	
	The latter assumption implies \begin{equation}\label{formula:image_of_sheaf}
		\widetilde{\varphi}^*_n S^*_{Y_n}= S^*_n
	\end{equation}
	for all $n$. We write equality, as such an isomorphism is determined up to a scalar $c_n$, and we assume that the scalars $c_n$ are chosen in a way compatible with the restriction
	maps
	\begin{equation*}
		\begin{tikzcd}
			\widetilde{\varphi}^*_n S^*_{Y_n}|_{X_{n-1}} \arrow[d,equal] \arrow{r}
			& \widetilde{\varphi}^*_{n-1}S^*_{Y_{n-1}} \arrow[d,equal] \\
			S^*_n|_{X_{n-1}} \arrow{r}
			& S^*_{n-1}.
		\end{tikzcd}
	\end{equation*}
	Then, by our above remark that isomorphisms of grassmannians are encoded by linear operators, the isomorphisms $\widetilde{\vfi}_n\colon X_n\stackrel{\sim}{\rightarrow} Y_n$ are recovered by a choice of compatible invertible linear operators \begin{equation*}
		\varphi^*_n\colon (V'_n)^*=\mathrm{H}^0(Y_n,S^*_{Y_n})\stackrel{\sim}{\rightarrow} \mathrm{H}^0(X_n,S^*_n)=V^*_n. 
	\end{equation*}
	The operators $\vfi^*_n$ are dual to unique operators $\vfi_n=(\vfi^*_n)^*\colon V_n \rightarrow V'_n$ which we will also consider.
	
	Since both chains of embeddings $X_n\hookrightarrow X_{n+1}$ and $Y_n\hookrightarrow Y_{n+1}$ are strict standard extensions, we have $\ilm V'_n=V=\ilm V_n$, and consequently, $\plm (V'_n)^*=V^*=\plm V^*_n$. Therefore, the operators $\varphi^*_n$ induce a linear operator \begin{equation*}
		\Phi\colon V^*\rightarrow V^*,
	\end{equation*}
	and more precisely, a commutative diagram
	
	
	
	\begin{equation}\label{formula:big_commutative_diagram}
		\begin{tikzcd}
			\begin{array}{c}
				V^*\\
				\vdots
			\end{array} \arrow[d]  
			& \begin{array}{c}
				V^*\\
				\vdots
			\end{array}   \arrow[d]  \arrow[yshift=2.0ex,l,"\Phi"'] \\
			V^*_n \arrow[d]
			& (V'_{n})^* \arrow[l,"\varphi^*_{n}"'] \arrow[d]\\
			V^*_{n-1} \arrow[d]
			& (V'_{n-1})^* \arrow[l,"\varphi^*_{n-1}"'] \arrow[d]\\
			\vdots & \vdots \\
		\end{tikzcd}
	\end{equation}

	The diagram (\ref{formula:big_commutative_diagram}) encodes the automorphism $\widetilde{\vfi}$ in the following way. Let $W'=\ilm(W'\cap V_n)$ be a point of $\Gra(W,E,V)$. Then \begin{equation}\label{formula:action_tilde_vfi_on_subspace}
		\widetilde{\vfi}(W')=\ilm \vfi_n(W'\cap V_n)=\ilm  \left(\left(\varphi_n^*\right)^{-1}\left((W'\cap V_n)^{\perp}\right)^{\perp}\right),
	\end{equation} where the orthogonal to $W'\cap V_n$ is taken in $V^*_n$ and the orthogonal to $\left(\varphi_n^*\right)^{-1}\left(\left(W'\cap V_n\right)^{\perp}\right)$ is taken in $V'_n$.
	
	Next, it is essential to observe that the subspace $V_{* E}\subset V^*$ is nothing but the subspace of global sections $\mu$ of the sheaf $\plm S^*_n$ satisfying the condition: the value of 
	$\mu$ at any point $W'\in \Gra(W,E,V)$ is a linear function on $W'$ which belongs to the subspace $W'_*\subset {W'}^*$. Here $W'_*$ is defined in terms of a basis of $W'$ which differs from $E$ by finitely many vectors ($W'_*$ is the span of the system of linear functionals dual to such a basis). Note that, for each $W'\in \Gra(W,E,V)$, the subspace $W'_*\subset W'^*$ is determined solely by the set $\Gra(W,E,V)$, and hence the above subspace of global sections $\mu$, i.e., the space $V_{*E}$, must be invariant under the operator $\Phi$. 
	Next, the subspace $V_E$ is the counterpart of the space $V_{*E}$ for the ind-grassmannian $\Gra(W^\perp,\wt E^*,V_*)$, where $V_*$ is defined by the fixed extension $\widetilde{E}$ of $E$ to a basis of $V$, and $E^*$ are the linear functions in $\widetilde{E}^*$ which do not vanish on $W$. 
	Since $\Gra(W^\perp,\wt E^*,V_*)$ is isomorphic to $\Gra(W,E,V)$, the space $V_E$ is also invariant under the linear map $\Phi'\colon (V_*)^*\rightarrow (V_*)^*$ induced by the automorphism $\widetilde{\varphi}$ of $\Gra(W,E,V)$. 
	
	We have shown that any automorphism $\widetilde{\varphi}\colon \Gra(W,E,V)\rightarrow \Gra(W,E,V)$ satisfying (\ref{formula:second_possibillity_for_varphi_tilde}) induces a pair of invertible operators $\varphi:=\Phi'|_{V_{E}}\colon V_E\rightarrow V_E$ and $\overline{\varphi}:=\Phi|_{V_{*E}}\colon V_{*E}\rightarrow V_{*E}$ which determine an element of the Mackey group $ G(V_E,V_{*E})$. Moreover, if $\mathrm{Aut}^0\Gra(W,E,V)$ stands for the group of automorphisms of the ind-variety $\Gra(W,E,V)$ satisfying (\ref{formula:second_possibillity_for_varphi_tilde}), then   the assignment $\widetilde{\varphi}\longmapsto \varphi$, or equivalently $\widetilde{\varphi}\longmapsto \overline{\varphi}^{-1}$, defines an injective group homomorphism $$\varepsilon \colon \mathrm{Aut}^0\Gra(W,E,V)\hookrightarrow \Pro G(V_E,V_{*E}).$$
	
	We now check that the action of the image in $\Pro G(V_E,V_{*E})$ of $\mathrm{Aut}^0\Gra(W,E,V)$ is given by the formula \begin{equation*}
		\widetilde{\varphi}\left(W'\right)=\overline{\varphi}^{-1}\left(W'^{\perp}\right)^{\perp}\cap V,
	\end{equation*}
	where $W'^{\perp}\subset V_{*E}$ and $\overline{\varphi}^{-1}\left( W'^{\perp}\right)^{\perp}\subset (V_{*E})^*$ (clearly, $V\subset (V_{*E})^*$). To do this, recall that  $\widetilde{\vfi}(W')$ is given by formula (\ref{formula:action_tilde_vfi_on_subspace}). 
	Therefore, we need to verify that \begin{equation}\label{formula:two_ilm_are_equal}
		\ilm \left(\left(\varphi^*_n\right)^{-1}\left(\left(W'\cap V_n\right)^\perp\right)^\perp\right)=\overline{\varphi}^{-1}\left(\left(\plm\left(\left(W'\cap V_n\right)^\perp\right)\right)^\perp\right) \cap V.
	\end{equation}
	However, formula (\ref{formula:two_ilm_are_equal})  follows from the observation that both its left-hand and right-hand sides coincide with the subspace of vectors in $V$ which vanish on $\left(\varphi^*_n\right)^{-1}\left(\left(W'\cap V_n\right)^\perp\right)$ whenever they belong to $V_n$; we consider here vectors in $V_n$ as 
	linear functions on $V^*_n$.
	
	In conclusion, the image of $\mathrm{Aut}^0\Gra(W,E,V)$ in $\Pro G(V_E,V_{*E})$ is a subgroup of $\Pro G(V_E,V_{*E})$ which acts on $\Gra(W,E,V)$ via the formula (\ref{formula:two_ilm_are_equal}). 
	Next, we note that since $\GL(\widetilde{E},V)$ acts transitively on $\Gra(W,E,V)$, for any $\widetilde{\vfi}\in \mathrm{Aut}^0\Gra(W,E,V)$ there are $\varkappa\in \GL(\widetilde{E},V)$ and  $\widetilde{\vfi}_W\in \mathrm{Aut}^0\Gra(W,E,V)$ 
	such that $\widetilde{\vfi}_W(W)=W$ and $\widetilde{\vfi}=\varkappa\circ \widetilde{\vfi}_W$. 
	Indeed if 
	$\widetilde{\vfi}(W)=W'$ for $W'\in \Gra(W,E,V)$, then $\widetilde{\vfi}=   \varkappa^{-1}\circ \varkappa \circ \widetilde{\vfi}$ where $\varkappa\in \GL(\widetilde{E},V)$ satisfies $\varkappa(W')=W$.
	Consequently,  $\vfi=\varkappa^{-1}\vfi_W$ for $\vfi_W\in \St_W$, in other words,   the image of $\varepsilon$ lies in $\Pro(\GL(\widetilde{E},V)\cdot \St_W)$. 
	
	To complete the proof, we need to show that any operator $\vfi \in \GL(\widetilde{E},V)\cdot \St_W$ determines a well-defined automorphism of $\Gra(W,E,V)$. In the Appendix we introduce the degree $d(\vfi)$ of an operator $\vfi$ and show that the space $\vfi \cdot W'=\overline{\vfi}^{-1}({W'}^{\perp})^{\perp}\cap V$ is $\widetilde{E}$-commensurable with $W$ for any $\vfi \in  G(V_E,V_{*E})$ such that $d(\vfi)=0$. Denote by $G^0(V_E,V_{*E})$ the group of all operators $\vfi\in G(V_E,V_{*E})$ with $d(\vfi)=0$.
	
	We have to convince ourselves that the action of $G^0(V_E,V_{*E})$
	on $\Gra(W,E,V)$ is by automorphisms of ind-varieties, and not merely by bijections of the set $\Gra(W,E,V)$. Let $V=W\oplus U$ where $\widetilde{E}\cap U$ spans $U$. Recall that our fixed nested finite-dimensional spaces $V_1 \subset \ldots \subset V_n \subset V_{n+1}\subset \ldots$ are spanned by elements of $\widetilde{E}$ and $\Gra(W,E,V)=\ilm \Gra(d_n,V_n)$ where $d_n=\dim V_n\cap W$. The embeddings $X_n=\Gra(d_n,V_n)\hookrightarrow \Gra(d_{n+1},V_{n+1})=X_{n+1}$ are strict standard extensions \begin{equation*}
		F_{d_n}\longmapsto F_{d_n}\oplus W_{n|n+1},
	\end{equation*}
	where $W_{n+1}=W\cap V_{n+1}=W_n\oplus W_{n|n+1}$ for $W_n=W\cap V_n$, and $W_{n|n+1}\cap \widetilde{E}$ spans $W_{n|n+1}$. Therefore we have a decomposition $V_{*E}={\overline{W}_{n}}_* \oplus V^*_n \oplus U^*_n$, where ${\overline{W}_{n}}_*$ is the $\widetilde{E}$-compatible direct complement of $V^*_n$ within $W_*+ V^*_n$ and $U^*_n$ is the dual of the $\widetilde{E}$-compatible direct complement $U_n$ of $W+V_n$ in $V$. Any invertible linear operator $\zeta\colon V_{*E}\rightarrow V_{*E}$ induces a decomposition\begin{equation*}
		V_{*E}=\zeta({\overline{W}_{n}}_*)\oplus \zeta(V^*_n)\oplus \zeta(U^*_n),
	\end{equation*} 
	and hence an operator\begin{equation*}
		\zeta_n\colon V^*_n:=V_{*E}/({\overline{W}_{n}}_*\oplus U^*_n)\longrightarrow V_{*E}/(\zeta({\overline{W}_{n}}_*)\oplus \zeta(U^*_n))=\zeta_n(V^*_n).
	\end{equation*}
	Moreover, we have $\zeta=\left(\plm\zeta_n\right)|_{V_{*E}}.$
	
	Each linear operator $\zeta_n$ induces an isomorphism of grassmannians \begin{equation*}
		\Gra(d_n,\zeta_n(V^*_n)^*)\longrightarrow\Gra(d_n,V_n),
	\end{equation*}
	and the varieties $\Gra(d_n,\zeta_n(V^*_n)^*)$ form an ind-variety isomorphic to $\Gra(W,E,V)$. We conclude that, if $\ilm \Gra(d_n,\zeta_n(V^*_n)^*)=\Gra(W,E,V)$ then $\zeta$ induces an automorphism of the ind-variety $\Gra(W,E,V)$. Set now $\zeta:=\overline{\vfi}^{-1}$ for $\vfi\in G^0(V_E,V_{*E})$. Then $\zeta_n=\left(\vfi^*_n\right)^{-1}$ and $\zeta_n(V^*_n)^*=\vfi_n(V_n)$.  By the above mentioned result from the Appendix, we know that  $\vfi \cdot W'\in \Gra(W,E,V)$ whenever $\vfi \in G^0(V_E,V_{*E})$ and $W'\in \Gra(W,E,V)$.
	Therefore $\ilm\Gra(d_n,\varphi_n(V_n))=\Gra(W,E,V)$,
	and we have shown that $\Pro G^0(V_E,V_{*E})\subset \Imm\varepsilon$. 
	
	Finally, Theorem A.\ref{theoapp:d_0_G_k_G_k} (i) from the Appendix implies that $\St_W\subset G^0(V_E,V_{*E})$, and hence that also $\GL(\widetilde{E},V) \cdot \St_W \subset G^0(V_E,V_{*E})$. Consequently, $\GL(\widetilde{E},V) \cdot \St_W=G^0(V_E, V_{*E})$ and $\GL(\widetilde{E},V) \cdot \St_W$ is a group. In particular, $\GL(\widetilde{E},V)\cdot \St_W=St_W\cdot \GL(\widetilde{E},V)$. The proof is now complete as we have shown that $\Imm\varepsilon=\Pro(\GL(\widetilde{E},V) \cdot\St_W ).$
	\hspace{\fill}$\square$\par
	\corop{If $\dim W=\codim_V W =\infty$, the group $\mathrm{Aut}\Gra(W,E,V)$ is isomorphic to the projectivization of the connected component of unity in the group Japanese $\GL(\infty)$.}{We proved that $\mathrm{Aut}\Gra(W,E,V)\cong \Pro G^0(V_E, V_{*E})$. As pointed out in the Appendix, the group $G^0(V_E, V_{*E})$ is isomorphic to the connected component of unity in the group Japanese $\GL(\infty)$.}
	

	\sect{Proof of Theorem \ref{theorem:1} in the general case}\label{sect:proof_theorem_1_in_the_general_case}
	
	\emph{Step 1. Reduction to the case of an automorphism which preserves all inverse limits of dual tautological bundles.} Consider the ind-variety $\Fl(W,\widetilde{E},V)$ for our fixed generalized flag $W=\{ W_{\alpha}\}$ compatible with the fixed basis $\widetilde{E}$ of $V$. Fix an exhaustion of $\Fl(W,\widetilde{E},V)$ as a direct limit $\ilm \Fl(\underline{d}_{n},V_{n})$ of strict standard extensions of finite-dimensional flag varieties $\Fl(\underline{d}_{n},V_{n})$. If the generalized flag $W$ is symmetric, then the exhaustion can be chosen so that each flag variety $\Fl(\underline{d}_{n},V_{n})$ is symmetric, i.e., the vector $(d_n^0=0,d^1_{n},d^2_{n},\ldots,d^s_{n},d_n^{s+1}=\dim V_{n})$ satisfies the condition $$d^i_{n}-d^{i-1}_{n}=d^{s+2-i}_{n}-d^{s+1-i}_{n}$$ for all $1\leq s\leq n+1$. If $W$ is not symmetric, then infinitely many vectors $\underline{d}_{n}$ are not symmetric, so (by passing to a subsequence of the sequence $\{n\}$) we can assume that all vectors $\underline{d}_{n}$ are not symmetric. 
	
	Let $\widetilde{\varphi}\colon \Fl(W,\widetilde{E},V) \rightarrow \Fl(W,\widetilde{E},V)$ be an automorphism. Set $X_n:=\Fl(\underline{d}_{n},V_{n})$ and let $Y_n=\widetilde{\varphi}(X_n)$. The varieties are finite-dimensional flag varieties, and there are fixed isomorphisms $\widetilde{\varphi}_n:=\widetilde{\varphi}|_{X_n}\colon X_n \stackrel{\sim}{\rightarrow} Y_n$. Moreover, Theorem \ref{theorem:P} implies that 
	the embeddings $Y_n\hookrightarrow Y_{n+1}$ are standard extensions. By replacing $Y_n$ or $Y_{n+1}$ by its dual flag variety we can further assume that the embeddings $Y_n\hookrightarrow Y_{n+1}$ are strict standard extensions. 
	
	Denote by $S^j_n$ and $S^j_{Y_n}$ the tautological bundles of rank $d^j_n$ on $X_n$ and $Y_n$, respectively. There are two possibilities: either
	\begin{equation}\label{formula:first_case_proof_theorem_1}
		\widetilde{\varphi}^*_n S^j_{Y_n} \simeq S^j_n
	\end{equation}
	for all $n$, or
	\begin{equation}\label{formula:second_case_proof_theorem_1}
		\widetilde{\varphi}^*_n S^j_{Y_n}  \simeq (\widetilde{V}_n/S^{s+1-j}_n)^*
	\end{equation}
	for some $n=n_0$ and all $j$, $1<j\lee s$. Case (\ref{formula:second_case_proof_theorem_1}) can occur only if the vector $\underline{d}_{n_0}$ is symmetric. Moreover, then (\ref{formula:second_case_proof_theorem_1}) will necessarily hold for all $n>n_0$ due to the assumption that all embeddings $X_{n}\hookrightarrow X_{n+1}$ and $Y_{n}\hookrightarrow Y_{n+1}$ are strict standard extensions. In that case, we may as well assume that (\ref{formula:second_case_proof_theorem_1}) holds for all $n$.
	
	Similarly to the case of an ind-grassmannian, if (\ref{formula:second_case_proof_theorem_1}) holds for all $n$ we can compose $\widetilde{\varphi}$ with an automorphism \begin{equation*}
		\delta: \Fl(W,\widetilde{E},V) \stackrel{\sim}{\longrightarrow} \Fl(W^{\perp},\widetilde{E}^*,V_*) \stackrel{\sim}{\longrightarrow} \Fl(W,\widetilde{E},V)
	\end{equation*}
	which maps $W$ first to $W^{\perp}$ and then maps $W^{\perp}$ to a point of $\Fl(W,\widetilde{E},V)$ under the isomorphism
	$ \Fl(W^{\perp},\widetilde{E}^*,V_*) \stackrel{\sim}{\longrightarrow} \Fl(W,\widetilde{E},V)$
	induced by an appropriate linear isomorphism $V_*\rightarrow V$ mapping $\widetilde{E}^*$ to $\widetilde{E}$. Then the composition $\delta\circ\widetilde{\varphi}$ satisfies the condition (\ref{formula:first_case_proof_theorem_1}) for all $n$. Therefore, in order to prove Theorem \ref{theorem:1}, it suffices to prove that \begin{equation*}
		\mathrm{Aut}^0\Fl(W,\widetilde{E},V)\simeq \Pro(\GL(\widetilde{E},V)\cdot \St_W)
	\end{equation*}
	where $\mathrm{Aut}^0\Fl(W,\widetilde{E},V)$ denotes the group of automorphisms of $\Fl(W,\widetilde{E},V)$ satisfying (\ref{formula:first_case_proof_theorem_1}) for all $n$.
	
	\emph{Step 2. From automorphisms to linear operators.} Since (\ref{formula:first_case_proof_theorem_1}) holds, the automorphism \begin{equation*}
		\widetilde{\varphi}\colon \Fl(W,\widetilde{E},V) \rightarrow \Fl(W,\widetilde{E},V)
	\end{equation*} induces automorphisms \begin{equation}\label{formula:induced_automorphisms}
		\widetilde{\varphi}_{\alpha}\colon \Gra(W_{\alpha},E_{\alpha},V)\rightarrow \Gra(W_{\alpha},E_{\alpha},V)
	\end{equation} 
	for each subspace $W_{\alpha}$ in $W$. In turn, the automorphisms (\ref{formula:induced_automorphisms}) induce linear operators, defined up to scalar multiples, \begin{equation*}
		\Phi_{\alpha}:V^*\rightarrow V^*
	\end{equation*}
	as explained in Section \ref{sect:proof_theorem_1_for_ind-grassmannians}.
	
	We now point out that the operators $\Phi_{\alpha}$ can be chosen to coincide, i.e., to define a single operator
	\begin{equation*}
		\Phi:V^*\rightarrow V^*
	\end{equation*}
	not depending on $\alpha$. This observation is justified as follows. Denote by $S^*_{\alpha}$ the pullback to $\Fl(W,\widetilde{E},V)$ of the inverse limit of the tautological bundles $S^*_{n}$ on $\Gra(W_{\alpha},E_{\alpha},V)$. Let $\alpha'<\alpha$ be two indices in the chain $W=\{W_{\alpha}\}$. Then the morphism of inverse limits $S^*_{\alpha}\rightarrow S^*_{\alpha'}$, arising from the respective morphism of inverse systems, induces a commutative diagram
	
	\begin{equation}\label{formual:diagram_with_H_and_S}
		\begin{tikzcd}
			\mathrm{H}^0(\Fl(W,\widetilde{E},V),S^*_{\alpha}) \arrow[d,equal] \arrow{r}
			& \mathrm{H}^0(\Fl(W,\widetilde{E},V),S^*_{\alpha'}) \arrow[d,equal] \\
			V^* \arrow[r,equal,"id"]
			& V^*,
		\end{tikzcd}
	\end{equation}
	the vertical equalities being the identifications $V^*=\plm V^*_n=\plm \mathrm{H}^0(\Gra(d^{\gamma}_n,V_n),S^*_n)$ where\linebreak $\Gra(W_{\gamma},E_{\gamma},V)=\ilm \Gra(d^{\gamma}_n,V_n)$ for $\gamma=\alpha$ and $\gamma=\alpha'$, respectively.
	Therefore, for a fixed $\alpha$ and all $\alpha'<\alpha$, the spaces of the form $\mathrm{H}^0(S^*_{\alpha'})$ are identified with $V^*$ in a way compatible with the upper horizontal arrows of the diagrams (\ref{formual:diagram_with_H_and_S}). Next, the following diagram is commutative
	\begin{equation*}
		\begin{tikzcd}
			\mathrm{H}^0(\Fl(W,\widetilde{E},V),S^*_{\alpha}) \arrow[d,equal] \arrow{r}
			\arrow[bend right = 70]{ddd} & \mathrm{H}^0(\Fl(W,\widetilde{E},V),S^*_{\alpha'})  \arrow[bend left=70]{ddd} \arrow[d,equal] \\
			V^* \arrow[r,equal,"id"] \arrow[d, "\Phi_{\alpha}"']
			& V^* \arrow[d,"\Phi_{\alpha'}"]\\
			V^* \arrow[r,equal,"id"] \arrow[d,equal]
			& V^* \arrow[d,equal] \\
			\mathrm{H}^0(\Fl(W,\widetilde{E},V),\widetilde{\varphi}^*S^*_{\alpha})  \arrow{r}
			& \mathrm{H}^0(\Fl(W,\widetilde{E},V),\widetilde{\varphi}^* S^*_{\alpha'}), 
		\end{tikzcd}
	\end{equation*}
	and this implies $\Phi_{\alpha}=\Phi_{\alpha'}$.
	
	\emph{Step 3. Injective homomorphism $\mathrm{Aut}^0 \Fl(W,\widetilde{E},V)\rightarrow \Pro G\left(V^W_{\widetilde{E}},V^W_{*\widetilde{E}}\right)$.} Note first that each homomorphism \begin{equation*}
		\mathrm{Aut}^0 \Fl(W,\widetilde{E},V) \longrightarrow \mathrm{Aut}^0 \Gra(W_{\alpha},\widetilde{E},V)
	\end{equation*}
	\begin{equation*}\widetilde{\varphi}\longmapsto \widetilde{\varphi}_{\alpha}
	\end{equation*}
	is injective, since $\widetilde{\varphi}_{\alpha}$ recovers $\widetilde{\varphi}$ through the formula \begin{equation*}
		\widetilde{\varphi}\left(W'\right)_{\alpha'}=\Phi^{-1}\left({W'_{\alpha'}}^{\perp} \right)^{\perp}\cap V
	\end{equation*}
	for any $\alpha'$, where ${W'_{\alpha'}}^{\perp}\subset V^{W_{\alpha'}}_{*E_{\alpha'}}$, $\Phi^{-1}\left({W'_{\alpha'}}^{\perp} \right)^{\perp} \subset V^{W_{\alpha'}}_{E_{\alpha'}}$ and, as explained above, $\Phi=\Phi_{\alpha}$. 
	
	Moreover, the automorphism $\Phi:V^*\rightarrow V^*$ induces pairs of automorphisms \begin{equation}
		\vfi_{\alpha}=\Phi'_{\alpha}|_{V^{W_{\alpha}}_{E_{\alpha}}}\colon V^{W_{\alpha}}_{E_{\alpha}}\longrightarrow\label{formula:phi_alpha} V^{W_{\alpha}}_{E_{\alpha}},~
		\overline{\vfi}_{\alpha}\colon V^{W_{\alpha}}_{*E_{\alpha}}\longrightarrow V^{W_{\alpha}}_{*E_{\alpha}},
	\end{equation}
	compatible with all inclusions of the form $V^{W_{\alpha'}}_{E_{\alpha'}}\subset V^{W_{\alpha}}_{E_{\alpha}}$, $V^{W_{\alpha}}_{*E_{\alpha}}\subset V^{W_{\alpha'}}_{*E_{\alpha'}}$  for $\alpha'<\alpha$, and such that $\vfi_{\alpha}
	\in  G\left(V^{W_{\alpha}}_{E_{\alpha}},V^{W_{\alpha}}_{*E_{\alpha}}\right)$. The compatibility of $\vfi_{\alpha}$ and $\vfi_{\alpha'}$ is clear, while the compatibility of $\overline{\vfi}_{\alpha}$ and $\overline{\vfi}_{\alpha'}$ follows from the inclusion
	\begin{equation*}
		\vfi_{\alpha'}^{-1}(W_{\alpha'}^{\perp})^{\perp}\subset\vfi_{\alpha}^{-1}(W_{\alpha}^{\perp})^{\perp},
	\end{equation*}
	where $W_{\alpha}^{\perp}\subset V^{W_{\alpha}}_{*E_{\alpha}}$, $W_{\alpha'}^{\perp}\subset V^{W_{\alpha'}}_{*E_{\alpha'}}$, $\vfi_{\alpha}^{-1}(W_{\alpha}^{\perp})^{\perp}\subset V^{W_{\alpha}}_{E_{\alpha}}$, $\vfi_{\alpha'}^{-1}(W_{\alpha'}^{\perp})^{\perp}\subset V^{W_{\alpha'}}_{E_{\alpha'}}$. Therefore, we conclude that the system of linear operators (\ref{formula:phi_alpha}) determines a unique element in $ G\left( V^W_{\widetilde{E}}, V^W_{*\widetilde{E}}\right)$
	for $V^W_{\widetilde{E}}=\bigcap \limits_{\alpha}V^{W_{\alpha}}_{E_{\alpha}}$ and $V^W_{*\widetilde{E}}=\bigcap \limits_{\alpha}V^{W_{\alpha}}_{*E_{\alpha}}$, and we obtain a homomorphism \begin{equation*}
		\varepsilon:\mathrm{Aut}^0 \Fl(W,\widetilde{E},V) \longrightarrow \Pro G\left(V^W_{\widetilde{E}},V^W_{*\widetilde{E}}\right).
	\end{equation*}
	The fact that $\GL(\widetilde{E},V)$ acts transitively on $\Fl(W,\widetilde{E},V)$ shows, by the same argument as in Section~\ref{sect:proof_theorem_1_for_ind-grassmannians}, that the image of $\varepsilon$ lies in the subset $\Pro\left(\GL(\widetilde{E},V)\cdot \St_W\right)$ of $\Pro G\left(V^W_{\widetilde{E}}, V^W_{*\widetilde{E}}\right)$.
	
	We have to prove that $\varepsilon$ is injective, i.e., that the image of $\varepsilon$ determines all operators $\vfi_{\alpha}
	$ and $\overline{\vfi}_{\alpha}$ 
	as above. We will do this by 
	recalling that each operator $\vfi_{\alpha}:V^{W_{\alpha}}_{E_{\alpha}}\rightarrow V^{W_{\alpha}}_{E_{\alpha}}$ admits a matrix as described in the Appendix. The key point is that if an invertible operator on $V^{W_{\alpha}}_{E_{\alpha}}$ admits such a matrix, then this matrix is unique. Now the compatibility of the operators $\vfi_{\alpha}$ under all inclusions $V^{W_{\alpha'}}_{E_{\alpha'}}\hookrightarrow V^{W_{\alpha}}_{E_{\alpha}}$ for $\alpha'<\alpha$ implies that the matrices of all operators $\vfi_{\alpha}$ coincide. Since each homomorphism $$\widetilde{\vfi}_{\alpha} \longmapsto \vfi_{\alpha}$$ is injective according to Section \ref{sect:proof_theorem_1_for_ind-grassmannians}, we conclude that $\varepsilon$ is injective.

	
	
	\emph{Step 4. The image of $\varepsilon$.} As a final step of the proof, we need to show that the image of $\varepsilon$ coincides with the set $\Pro(\GL(\widetilde{E},V)\cdot \St_W)$. For this it suffices to prove that $\Pro(\GL(\wt E,V)\cdot \St_W)$ 
	belongs to the image of $\varepsilon$, i.e., that $\GL(\wt E,V)\cdot\St_W$ acts on the ind-variety $\Fl(W,\widetilde{E},V)$ via the formula (\ref{formula:action_on_flag}). 
	
	
	Pick an operator $\vfi\in \GL(\wt E,V)\cdot\St_W$ and a flag $W'\in\Fl(W,\widetilde{E},V)$. Since $W$ and $W'$ are $\widetilde E$-commensurable, one has $\GL(\wt E,V)\cdot \St_{W}=\GL(\widetilde{E},V)\cdot \St_{W'}$. 
	Therefore $\vfi=\varkappa^{-1} \vfi_{W'}$ for some $\vfi_{W'}\in \St_{W'}$  and some $\varkappa\in \GL(\widetilde{E},V)$. Consequently,   $\vfi(W')=\varkappa^{-1}(W')$, 
	i.e. $\vfi^{-1}(W')$ is $\wt E$-commensurable with $W$. In conclusion, $\mathrm{im} \varepsilon=\Pro(\GL(\wt E,V)\cdot \St_{W})$, and since $\GL(\wt E,V)\cdot \St_{W}$ is a group we have also $\mathrm{im} \varepsilon=\Pro( \St_{W} \cdot\GL(\wt E,V))$. The proof is complete.
	
	
	\sect{An explicit matrix form of the group $\GL(\widetilde{E},V)\cdot \St_W $}\label{sect_explicit_matrix_form}
	
	Now we would like to characterize the product $\GL(\widetilde{E},V)\cdot \St_W$ in terms of matrices. We start by describing a matrix form of the group $\St_W$.
	
	Choose a linear order on $\widetilde{E}$ such that $e_j\in W_{\alpha} \setminus W_{\alpha'}$, $e_k\in W_{\alpha'}$ for $\alpha'<\alpha$ implies $k<j$.  It follows from the Appendix that for each space $W_{\alpha}$ the stabilizer of $W_{\alpha}$ in $ G \left( V^{W_{\alpha}}_{E_{\alpha}}, V^{W_{\alpha}}_{*E_{\alpha}}\right)$ under the action (\ref{formula:action_on_flag}) can be represented by infinite matrices (with rows and columns ordered by the fixed order on $\widetilde{E}$) 
	which, together with their inverses, have the form \begin{equation}\label{formula:matrix_form_zero_up_right}
		\left(\begin{array}{c|c}
			A & B  \\
			\hline
			0 & D  \end{array}\right)
	\end{equation}
	where $A$ has finitary rows, $D$ has finitary columns, and there are no restrictions on the rows and columns of $B$. Certainly, the splitting (\ref{formula:matrix_form_zero_up_right}) depends on the space $W_{\alpha}$.

	Next, the fact that the operators $\varphi_{\alpha}$ and $\overline{\varphi}_{\alpha}$ form a system compatible with the inclusions $V^{W_{\alpha'}}_{E_{\alpha}}\subset V^{W_{\alpha}}_{E_{\alpha}}$, $V^{W_{\alpha}}_{*E_{\alpha}}\subset V^{W_{\alpha'}}_{*E_{\alpha}}$ for $\alpha'<\alpha$, 
	implies that all operators $\vfi_{\alpha}$, respectively $\overline{\vfi}_{\alpha}$, are represented by the same matrix. 
	This means that all matrices (\ref{formula:matrix_form_zero_up_right}) are just one matrix  which satisfies the above conditions for all spaces $W_{\alpha}$. Consequently, 
	$\St_W$ consists of matrices $M$ which, together with their inverses, have the form
	\begin{equation}\predisplaypenalty=0\label{formula:matrix_form}
		M=\begin{pmatrix} 
			\ddots & & &  &   \\
			0 & \boxed{A_{\alpha'}} & &  &   \\
			&&\ddots & & \\
			0 & 0 &  & \boxed{A_{\alpha}} &   &  \\ 0
			& 0 &  & 0  & \ddots   \\
		\end{pmatrix},
	\end{equation}
	the rows and columns of $M$ being ordered by the ordered set which orders the elements of $\widetilde{E}$ as above, 
	and the diagonal blocks being of size $\dim (W_{\alpha_2}/W_{\alpha_1}) \times\dim (W_{\alpha_2}/W_{\alpha_1})$ where $\alpha$ is one of the indices $\alpha_1$ and $\alpha_2$, and $\alpha_1$ is the immediate predecessor of $\alpha_2$. 
	This datum induces a block structure on the entire matrix $M$, and all strictly lower-triangular blocks are set to equal zero.
	Furthermore, since $M$ and $M^{-1}$ are subject to the above 
	additional conditions for all splittings arising from spaces $W_{\alpha}$ of $W$, the matrices $M$ and $M^{-1}$ satisfy the following: 
	\begin{itemize}
		\item if there exists a rightmost highest  block $L$ (this depends on the 
		order on $\widetilde{E}$, i.e., ultimately on the structure of the generalized flag $W$), then there are no conditions on the rows and columns of $L$, all columns of $M$ which intersect  $L$ have at most finitely many nonzero entries outside of $L$, all rows of $L$ which intersect $L$ have at most finitely many nonzero entries outside of $L$, and all other rows and columns of $M$ have at most finitely many nonzero entries in the region 
		in which the first index is greater then the first index of any entry of $L$ and the second index is smaller than the second index of any entry of $L$;
		\item if $M$ has no rightmost highest block, then all rows 
		have at most finitely many nonzero entries in  direction to the left (from any point on) and all columns 
		have at most finitely many nonzero entries in the downward direction  (from any point on).
	\end{itemize}
	
	In order to pass to the full group $\GL(\widetilde{E},V) \cdot \St_W$, consider matrices of the form (\ref{formula:matrix_form}) and replace the zeros in the lower-triangular part by finitely many nonzero entries. Let's refer to such matrices as $W$-\emph{aligned}. If $M$ is a $W$-aligned matrix, then every space $W_\alpha$ induces a splitting of $M$ into four blocks \begin{equation}\label{formula:splitting_depends_on_alpha}
		\left(\begin{array}{c|c}
			A & B  \\
			\hline
			C & D  \end{array}\right)
	\end{equation}
	where $C$ is a finitary matrix. A $W$-aligned matrix M is \emph{eligible} if $M^{-1}$ is also $W$-aligned, and for any $\alpha$ the splittings 
	\begin{equation*}
		\left(\begin{array}{c|c}
			A & B  \\
			\hline
			C & D  \end{array}\right)\text{ and }     \left(\begin{array}{c|c}
			A' & B'  \\
			\hline
			C' & D'  \end{array}\right)
	\end{equation*}
	of $M$ and $M^{-1}$ respectively, satisfy the condition \begin{equation*}
		\mathrm{rk}\,C=\mathrm{rk}\,C'.
	\end{equation*} 
	We leave it as an exercise to the reader to check that the condition of eligibility is empty (i.e., it is automatically satisfied) if the ordered set parameterizing the rows and columns of $A$ or $D$ is finite.

	Our result in this section states as follows.
	
	\theop{\label{theorem:description_general_case}The group $\GL(\widetilde{E},V) \cdot \St_W$ is isomorphic to the group of all eligible $W$-aligned matrices.}{Any matrix in $\GL(\widetilde{E},V)\cdot \St_W$ has the form $M_f M$ where $M_f$ is an element of $\GL(\widetilde{E},V)$ and $M$ is a matrix from $\St_W$ having the form (\ref{formula:matrix_form}). The necessary and sufficient condition for $M_f M$ to lie in $\GL(\widetilde{E},V)\cdot \St_W$ is that the splitting (\ref{formula:splitting_depends_on_alpha}) of $M_f M$ for any space $W_{\alpha}$ in $W$ satisfies the condition of Theorem A.\ref{theoapp:d_0_G_k_G_k}. For all $\alpha$ these conditions amount precisely to the requirement that the matrix $M_f M$ be eligible.}

	\sect{The isotropic case}\label{sect:proof_theorem_2_isotropic_case}
	\theopwop{\label{theorem:2} Let $X=\FlO(W,\widetilde{E},V)$ or $\FlS(W,\widetilde{E},V)$ for some isotropic generalized flag $W$ compatible with an isotropic basis $\widetilde{E}$ of $V$. 
		We assume in addition that $W$ 
		does not have the form $0\subset W_1 \subset W^{\perp}_1 \subset V$ where $\dim W_1=1$ in the symplectic case\textup, and that $\dim W^{\perp}_1/W_1\neq 2$ in the orthogonal case if the basis $\widetilde{E}$ has the form \textup{(\ref{formula:second_basis})}.  Then \begin{equation*}
			\mathrm{Aut}X\simeq \Pro(\Orth(\widetilde{E},V)\cdot \St^{\Orth}_W)
		\end{equation*}or\begin{equation*}
			\mathrm{Aut}X\simeq \Pro(\Sp(\widetilde{E},V) \cdot \St^{\Sp}_W),
		\end{equation*} where $\St^{\Orth}_W$ and $\St^{\Sp}_W$ denote respectively the stabilizer of $W$ in $\Orth(V)$ and $\Sp(V)$. 
		The action of $\Orth(V)$ on\linebreak $\FlO(W,\widetilde{E},V)$\textup, or respectively of $\Sp(V)$ on $\FlS(W,\widetilde{E},V)$\textup, is induced by the linear action of $\Orth(V)$ or $\Sp(V)$ on $V$.}

	\par\textsc{Proof}. 
	\emph{Step 1. The case of an isotropic ind-grassmannian.} The argument starts in the same way as for ordinary ind-grassmannians. The embeddings $X_n\hookrightarrow X_{n+1}$ are assumed to be standard extensions of isotropic grassmannians, and by Theorem \ref{theorem:P} 
	the embeddings  $Y_{n}:=\varphi(X_n)\xhookrightarrow{\eta_n} Y_{n+1}:=\varphi(X_{n+1})$ are also standard extensions. 
	
	Next, in the isotropic case the isomorphisms (\ref{formula:image_of_sheaf}) must hold, so we arrive to diagram (\ref{formula:big_commutative_diagram}) or, equivalently, to a commutative diagram 
	
	\begin{equation*}
		\begin{tikzcd}
			\begin{array}{c}
				V=\ilm V_n \\
				\vdots
			\end{array} \arrow[r,yshift=2.0ex,"\ilm\vfi_n"] 
			& 
			\begin{array}{c}
				V=\ilm V'_n  \\
				\vdots
			\end{array}     \\
			V_n \arrow[u] \arrow[r,"\varphi_{n}"]
			& V'_{n}  \arrow[u]\\
			V_{n-1} \arrow[u] \arrow[r,"\varphi_{n-1}"]
			& V'_{n-1}  \arrow[u]\\
			\vdots \arrow[u] & \vdots \arrow[u] \\
		\end{tikzcd}.
	\end{equation*}
	The isomorphism $V=\ilm V'_n$ holds since the embeddings $Y_n \hookrightarrow Y_{n+1}$ are standard extensions, and hence any $V'_n$ is a subspace of $V_N$ for some $N>n$. Furthermore, the restriction of the fixed (symmetric or antisymmetric) form on $V$ to $V_n'$ is a form defining the corresponding isotropic grassmannian in $V_n'$. (Such a form is unique up to a scalar if $W$ satisfies the conditions of the theorem.) Recall (from Section \ref{sect:proof_theorem_1_for_ind-grassmannians}) that the maps $\vfi_n$ are defined up to compatible scalars $c_n$.
	The key observation is that there is a unique choice of these scalars such that the maps $\vfi_n$ are isomorphisms of orthogonal or, respectively, symplectic vector spaces. This follows from the fact that the automorphism groups of our finite-dimensional isotropic grassmannians $X_n$ and $Y_n$ are the respective orthogonal or  symplectic groups. 
	
	Set  
	\begin{equation*}
		\vfi:=\ilm\vfi_n\colon V\rightarrow V.
	\end{equation*}
	Then by construction $\vfi$
	is an orthogonal or, respectively, symplectic operator, and $\vfi$ determines our automorphism $\widetilde{\varphi}$ which acts on a point $W'$ by formula (\ref{formula:action_tilde_vfi_on_subspace}). However, in the case considered we have \begin{equation}\label{formula:phi^-1_V_as_a_inductive_limit}
		\ilm \left(\varphi_n\left(\left(W'\cap V_n\right)^\perp\right)^\perp\right)=\ilm \left(\varphi_n\left(W'\cap V_n\right)\right)=\varphi(W'),
	\end{equation}
	hence the action of $\widetilde{\varphi}$ on $X$ is simply induced by the action of $\varphi$ as a linear automorphism of $V$. This implies that there is an injective homomorphism \begin{equation*}
		\varepsilon:\mathrm{Aut}\GraO(W,\widetilde{E},V)\hookrightarrow \Pro\Orth(V)
	\end{equation*} or, respectively,
	\begin{equation*}
		\varepsilon:\mathrm{Aut}\GraS(W,\widetilde{E},V)\hookrightarrow \Pro\Sp(V), 
	\end{equation*}
	\begin{equation*}\widetilde{\varphi}\longmapsto \varphi
	\end{equation*}
	and that the action of the image of $\varepsilon$ on $\GraO(W,\widetilde{E},V)$ or $\GraS(W,\widetilde{E},V)$ is induced by the linear action of $\Orth(V)$ or, respectively, $\Sp(V)$ on $V$.
	
	In the case of a general ind-variety of isotropic generalized flags, an injective homomorphism~$\varepsilon$ is constructed exactly as in Step 3 of the proof of Theorem \ref{theorem:1}. 
	Here the images of all homomorphisms $\varepsilon_{\alpha}$ lie in $\Pro \Orth(V)$ or, respectively, $\Pro \Sp(V)$, therefore $\varepsilon$ is just the homomorphism into the intersection of all images of $\varepsilon_{\alpha}$. 
	Moreover, the image of $\varepsilon$ coincides respectively with $\Pro(\Orth(\widetilde{E},V)\cdot \St^{\Orth}_W)$ or $\Pro(\Sp(\widetilde{E},V)\cdot \St^{\Sp}_W)$. Indeed, since $\Orth(\widetilde{E},V)$ or, respectively, $\Sp(\widetilde{E},V)$, acts transitively on $\FlO(W,\widetilde{E},V)$ or $\FlS(W,\widetilde{E},V)$, the image of $\varepsilon$ must be a subgroup of $\Pro(\Orth(\widetilde{E},V)\cdot \St^{\Orth}_W)$ or, respectively, $\Pro(\Sp(\widetilde{E},V)\cdot \St^{\Sp}_W)$. On the other hand, by the same argument as in Section \ref{sect:proof_theorem_1_in_the_general_case}, $\Orth(\widetilde{E},V)\cdot \St^{\Orth}_W$ or, respectively, $\Sp(\widetilde{E},V)\cdot \St^{\Sp}_W$ acts on the ind-variety $\FlO(W,\widetilde{E},V)$ or $\FlS(W,\widetilde{E},V)$ via the formula (\ref{formula:action_tilde_vfi_on_subspace}). Since in our case the equality (\ref{formula:phi^-1_V_as_a_inductive_limit}) holds, the proof is complete.\hspace{\fill}$\square$\par
	
	\corop{\label{corollary:group_structure_isotropic_case}The group $\Orth(\widetilde{E},V)\cdot \St^{\Orth}_W$, or $\Sp(\widetilde{E},V)\cdot \St^{\Sp}_W$, is isomorphic to the group of all invertible $W$-aligned matrices $M$ with finitary rows and columns satisfying $\overline{M}=\pm M^{-1}$, where $\overline{\,\,\cdot\,\,}$ denotes reflection along the antidiagonal, plus corresponds to the case of $O$ and minus corresponds to the case of $Sp$. 
	}{
		The group $\Orth(\widetilde{E},V)\cdot \St^{\Orth}_W$ or $\Sp(\widetilde{E},V)\cdot \St^{\Sp}_W$ is clearly the intersection of the group $\GL(\widetilde{E},V)\cdot \St_W$ with $\Orth(V)$ or, respectively, with $\Sp(V)$, and this implies the claim. Note that the condition $\overline{M}= \pm M^{-1}$ makes $M$  eligible automatically.}
	

	\begin{center}
		\subsection*{Appendix: on the structure of the group $ G (V_E,V_{*E})$}
	\end{center}
	
	In this Appendix we have collected some basic statements about the Mackey group $ G (V_E,V_{*E})$. The isomorphism of this group a the group sometimes referred to as Japanese $GL(\infty)$ is stated at the end.
	
	\setcounter{snomer}{0}
	\renewcommand{\thesnomer}{\arabic{snomer}}
	
	\sstappendix{ A lemma from linear algebra}
	\lemmpappendix{Let $A$ be\label{lemm:fin_dim_Im} a countable-dimensional vector space\textup, and $\mu\colon A^*\to A$\textup, $\nu\colon A^*\to A$ be linear maps such that $\alpha(\mu(\beta))=\beta(\nu(\alpha))$ for all $\alpha,~\beta\in A^*$. Then the dimensions of the images of $\mu$ and $\nu$ are finite and equal.
	}{We first check that 
		$\dim\Imm\mu<\infty$. This will imply that $\dim\Imm\nu$ is also finite because the conditions on $\mu$ and $\nu$ are symmetric. Pick a basis $\{f_1,~f_2,~\ldots\}$ of $A$, and let $\{f_1^*,~f_2^*,~\ldots\}$ be the dual system in $A^*$. Assume that $\Imm\mu$ is infinite dimensional. Then there exist linear functions $\lambda_n\in A^*$ for $n\in\Zp_{>0}$, such that $\dim\langle\mu(\lambda_1),~\ldots,~\mu(\lambda_n)\rangle_{\Cp}=n$. It is easy to see that these linear functions can be chosen so that $\lambda_n(f_i)=0$ for $i<n$ and $\lambda_n(f_n)\neq0$. Let $B:=\langle \lambda_n | n\in\Zp_{>0} \rangle_{\Cp}$. 
		The formula $(a,b)\mapsto b(a)$ for $a\in A$, $b\in B$ defines a non-degenerate pairing between $A$ and $B$. By a result of G. Mackey \cite[Lemma on p. 171]{Mackey45}, we may assume without loss of generality that $\lambda_n=f_n^*$ for all $n$.
		
		Since $\dim \Imm \mu = \infty$, for any $n\in\Zp_{>0}$ there exists $i_n\in\Zp_{>0}$ such that $\mu(f_{i_n}^*)$ is not contained in~$A_n:=\langle f_1,\ldots,f_n\rangle_{\Cp}$.
		In other words, there exist two infinite sequences of integers $1\leq i_1<i_2<\ldots$ and $2\leq k_1<k_2<\ldots$ so that $\mu(f_{i_n}^*)\in A_{k_n}\setminus A_{k_n-1}$ for all $n\geq1$.
		
		Now we define a sequence $c_1,~c_2,~\ldots$ of complex numbers inductively as follows. 
		Set $c_1$ to be an arbitrary nonzero scalar. For $n>1$ we let $c_n$ to be an arbitrary complex number such that $$\sum_{j=1}^nc_jf_{k_j}^*(\mu(f_{i_n}^*))=c_nf_{k_n}^*(\mu(f_{i_n}^*))+\sum_{j=1}^{n-1}c_jf_{k_j}^*(\mu(f_{i_n}^*))\neq0.$$ (Such a number exists because $f_{k_n}^*(\mu(f_{i_n}^*))\neq0$.) Finally, we define $\alpha\in A^*$ by the formula
		\begin{equation*}
			\alpha(f_j)=\begin{cases}c_n&\text{if }j=k_n\text{ for some }n\geq1,\\
				0&\text{otherwise}.
			\end{cases}
		\end{equation*}
		Then 
		\begin{equation*}
			\begin{split}f_{i_n}^*(\nu(\alpha))=\alpha(\mu(f_{i_n}^*))=\sum_{j=1}^nc_jf_{k_j}^*(\mu(f_{i_n}^*))\neq0
			\end{split}
		\end{equation*}
		for any $n\in \Zp_{>0}$. This means that $\nu(\alpha)\notin A_{i_n-1}$ for all $n\geq1$, i.e., that $\nu(\alpha)$ can not be expressed as a finite linear combination of $f_1,~f_2,~\ldots$, and this is a contradiction.
		
		Hence, $\dim\Imm \mu<\infty$, and also $\dim\Imm \nu<\infty$. It remains to check that $\dim\Imm \mu=\dim \Imm \nu$. 
		Since $\dim \Imm \mu<\infty$, there exist $n\in\Zp_{>0}$ and linear functions $\alpha_1,~\ldots,~\alpha_n\in A^{**}$ such that $\mu(a)=\sum_{i=1}^n\alpha_i(a)f_i$ for all $a\in A^*$. 
		We note that in fact $\alpha_i\in A$ for all $i$. Indeed, $$f_j^*(\mu(a))=f_j^*\left(\sum_{i=1}^n\alpha_i(a)f_i\right)=\alpha_j(a)=a(\nu(f_j^*)),$$ so $\alpha_j=\nu(f_j^*)\in A$ for all $j$. It follows that there exist $m\in\Zp_{>0}$ and scalars $\alpha_{ij}\in\Cp$ so that $\alpha_j=\sum_{i=1}^m\alpha_{ij}f_i$ for all $1\leq j\leq n$. Clearly, $\dim\Imm\mu=\rk(\alpha_{ij})_{i,j=1}^{m,n}$. Similarly, there exist $n',~m'\in\Zp_{>0}$, vectors $\beta_j\in A$, and scalars $\beta_{ij}\in\Cp$ so that $\nu(b)=\sum_{j=1}^{n'}\beta_j(b)f_j$ and $\beta_j=\sum_{i=1}^{m'}\beta_{ij}f_i$ for all $1\leq j\leq n'$. Then $\dim\Imm\nu=\rk(\beta_{ij})_{i,j=1}^{m',n'}$. However, $\alpha_{ij}=f_i^*(\alpha_j)=f_i^*(\nu(f_j^*))=f_j^*(\mu(f_i^*))=f_j^*(\beta_i)=\beta_{ji}$ for all $i,~j$. This means that $m'=n$, $n'=m$, and the matrices $(\alpha_{ij})_{i,j=1}^{m,n}$ and $(\beta_{ij})_{i,j=1}^{n,m}$ are transpose to each other. The result follows.}
	
	\sstappendix{A grading on the group $ G (V_E, V_{*E})$} We work in the setup of Section \ref{sect:proof_theorem_1_for_ind-grassmannians}. 
	Recall that we consider a generalized flag of the form $\{0\}\subset W\subset V$, where $\dim W=\codim_V W=\infty$, and that $E=\widetilde{E}\cap W$ is a basis of $W$. Let $U=\langle\widetilde{E}\setminus E\rangle_{\Cp}$. We have $V_{E}=(W_*)^*\oplus U$, $V_{*E}=W_*\oplus U^*$ and there is a natural nondegenerate pairing $V_{E}\times V_{*E}\rightarrow \mathbb{C}$. The group $ G (V_E,V_{*E})$ consists of invertible linear operators $\vfi \colon V_E\rightarrow V_E$ such that $\overline{\vfi}=\restr{\vfi^*}{V_{*E}}\colon V_{*E} \rightarrow V_{*E}$ is a well-defined isomorphism. 
	
	Choose an ordering of $E$ via $\mathbb{Z}_{<0}$ and an ordering of $\widetilde{E}\setminus E$ via $\mathbb{Z}_{>0}$. Set
	\begin{equation*}
		E^k:=\begin{cases} E\setminus \{e_k,\ldots,e_{-1}\}\text{ for }k<0 \\ E\cup \{e_1,\ldots,e_k\}\text{ for }k>0  \end{cases}
	\end{equation*}
	and $W^k:=\langle E^k \rangle_{\Cp}$. The disjoint union $\bigsqcup\limits_{k\in \Zp}\Gra(W^k,E^k,V)$, where $E=E^0$, $W=W^0$, is an ind-variety which we denote by $X$.

	\proppappendix{The formula \textup{(\ref{formula:action_on_flag})} defines an action of $ G (V_E, V_{*E})$ on $X$ by ind-variety\label{propapp:action_X} automorphisms.}{ We first prove that $\varphi\cdot W'\in X$ for any $W'\in X$ and any $\varphi\colon V_{E}\rightarrow V_{E}$, $\varphi\in  G (V_E, V_{*E})$. Since 
		in the definition of $X$, $W$ can be replaced by $W'$, it suffices to check that $\varphi\cdot W\in X$. 
		Denote by $\pi_U\colon V_E\to U$ and $\pi_{W_*}\colon V_{*E}\to W_*$ the canonical projections onto the corresponding direct summands. Set also $\mu:=\pi_U\circ \left(\vfi|_{(W_*)^*}\right)$ and $\nu:=\pi_{W_*}\circ (\overline{\vfi}|_{U^*})$. 
		
		We claim that $\dim\Imm\mu$ and $\dim\Imm\nu$ are finite (and equal). 
		Fix an isomorphism of vector spaces $\eta\colon W_*\to U$  
		and the dual isomorphism $\eta^*\colon U^*\to (W_*)^*$. 
		Consider the linear maps $\mu'=\mu\circ\eta^*$ and $\nu'=\eta\circ\nu$ from $U^*$ to $U$. Then, for all $\alpha,~\beta\in U^*$, we have
		\begin{equation*}
			\begin{split}\alpha(\mu'(\beta))&=\alpha(\mu(\eta^*(\beta)))=\alpha(\pi_U(\vfi(\eta^*(\beta))))=\alpha(\vfi(\eta^*(\beta)))=
				\overline{\vfi}(\alpha)(\eta^*(\beta))\\
				&=\pi_{W_*}(\overline{\vfi}(\alpha))(\eta^*(\beta))=\nu(\alpha)(\eta^*(\beta))=\eta(\nu(\alpha))(\beta)=\beta(\nu'(\alpha)).
			\end{split}
		\end{equation*}
		Hence, $\dim\Imm\mu'=\dim\Imm\nu'<\infty$ (and, consequently, $\dim\Imm\mu=\dim\Imm\nu<\infty$) by Lemma~A.\ref{lemm:fin_dim_Im}.
		
		Now, denote by $A_0$  the preimage of zero under the restriction of $\pi_U$ to $\vfi((W_*)^*)$. Obviously, $A_0=(W_*)^*\cap\vfi((W_*)^*)$. According to the proof of Lemma~A.\ref{lemm:fin_dim_Im}, there exist $n\in\Zp_{>0}$ and $\alpha_1,~\ldots,~\alpha_n\in W_*$ such that $$\mu(\omega)=\pi_U(\vfi|_{(W_*)^*}(\omega))=\sum\nolimits_{i=1}^n\omega(\alpha_i)e_i$$ for all $\omega\in (W_*)^*$. A vector $\omega$ of $(W_*)^*$ belongs to $A_0$ if and only if $\mu(\omega)=0$, i.e, if $\omega(\alpha_i)=0$ for $1\leq i\leq n$. 
		It follows immediately that $A_0\cap W$ has finite codimension in $W$ and, moreover, given $w\in W$, one has $w\in A_0$ if and only if $\alpha_i(w)=0$ for all $1\leq i\leq n$. Hence, $A_0\cap W$ contains all but finitely many vectors of $E$.
		
		Recall that $\vfi\cdot W=\overline{\vfi}^{-1}(W^{\perp})^{\perp}\cap V$. Clearly, $W^{\perp}=U^*$, and since $\overline{\vfi}^{-1}(B)^{\perp}=\vfi(B^{\perp})$ for any subspace~$B$ of $V_{*E}$, we conclude that $$\overline{\vfi}^{-1}(U^*)^{\perp}=\vfi((U^*)^\perp)=\vfi((W_*)^*).$$ Hence $\vfi\cdot W=\vfi((W_*)^*)\cap V$. Since $A_0\subset\vfi((W_*)^*)$, we have $A_0\cap W\subset\vfi\cdot W$, and so $\vfi\cdot W$ contains all but finitely many vectors of $E$. On the other hand, $$\vfi((W_*)^*)\subset (W_*)^*\oplus\pi_U(\vfi((W_*)^*))=(W_*)^*\oplus\Imm\mu.$$ This implies $\vfi\cdot W\subset W\oplus\Imm\mu$, 
		and consequently $\vfi\cdot W\in X$.
		
		Next, a trivial modification of the argument at the end of the proof of Theorem \ref{theorem:1} for ind-grassmannians shows that $\vfi$ induces an automorphism of ind-varieties\begin{equation*}
			\widetilde{\vfi}:X\rightarrow X.
		\end{equation*}
		Finally, the fact that $\widetilde{\vfi'\circ \vfi''}=\widetilde{\vfi'}\circ \widetilde{\vfi}''$ follows from the realization of $\tilde{\vfi}$ as an inverse limit of morphisms of finite-dimensional grassmannians where the property $\widetilde{\vfi_n'\circ \vfi_n''}=\widetilde{\vfi}_n'\circ \widetilde{\vfi}_n''$ is obvious.
	}

	The ind-variety $X$ is the disjoint union of the ind-grassmannians $\Gra(W^k,E^k,V)$ for $k\in\Zp$. These closed ind-subvarieties are not stable with respect to this action. 
	
	\examappendix{Let a \emph{reverse sequence} be a set parameterized by $\mathbb{Z}_{<0}$. The vectors in $V_E$ and $V_{*E}$ can be written respectively as $x=\sum \limits_{i>0} \omega_i e_{-i}+\sum \limits_{i>0} u_i e_{i}\in V_E$, $y=\sum \limits_{i>0} w_i e^*_{-i}+\sum \limits_{i>0} \nu_i e^*_{i}\in V_{*E}$  where $\omega=(\ldots,~\omega_{-2},~\omega_{-1})$ and $w=(\ldots,~w_{-2},~w_{-1})$ are reverse sequences of complex numbers, and $u=(u_1,~u_2,~\ldots)$, $\nu=(\nu_1,~\nu_2,~\ldots)$  are usual sequences of complex numbers. In addition, $u$ and $w$ are finitary, i.e. have at most finitely many nonzero entries. In what follows, we write simply $(\omega,u)$ for vectors in $V_E$, and $(w,\nu)$ for vectors in $V_{*E}$.
		
		Consider the shift (linear) operator $Sh\colon V_{E}\to V_{E}$ 
		\begin{equation*}
			Sh((\omega,~u))=((\ldots,~\omega_{-3},~\omega_{-2}),(\omega_{-1},~u_1,~u_2\ldots)).
		\end{equation*}
		One checks immediately that $Sh$ is an element of the group $ G (V_E,V_{*E})$ with $\overline{Sh}$ having the form $$\overline{Sh}((w,~\nu))=((\ldots,~w_{-2},~w_{-1},\nu_1),(\nu_2,~\nu_3,~\ldots)).$$
		Moreover, $Sh\cdot W$ belongs to $\Gra(W^1,E^1,V)$, and in fact
		$Sh\cdot\Gra(W,E,V)=\Gra(W^1,E^1,V)$.
		In addition, we note that
		\begin{equation*}
			\dim\pi_U(Sh((W_*)^*))-\dim\pi_U(Sh^{-1}((W_*)^*))=\dim\pi_{W_*}(\overline{Sh}(U^*))-\dim\pi_{W_*}(\overline{Sh}^{-1}(U^*))=1.
		\end{equation*}
	}
	
	This example motivates the following.
	
	\defiappendix{
		Let $\vfi\in  G (V_E, V_{*E})$. We\label{defi:defect} define the integer
		\begin{equation*}
			d(\vfi)=d(\overline{\vfi}):=\dim\pi_U(\vfi((W_*)^*))-\dim\pi_U(\vfi^{-1}((W_*)^*))=\dim\pi_{W_*}(\overline{\vfi}(U^*))-\dim\pi_{W_*}(\overline{\vfi}^{-1}(U^*))
		\end{equation*}
		to be the \emph{degree} of $\vfi$. (The latter equality follows immediately from the fact that $\dim\Imm(\pi_U\circ\restr{\vfi}{(W_*)^*})=\dim\Imm(\pi_{W_*}\circ\restr{\overline{\vfi}}{U^*})$ checked in the proof of Proposition A.\ref{propapp:action_X}.) 
	}
	
	\proppappendix{Given $\vfi\in  G (V_E,V_{*E})$ and $k\in\Zp$, one \label{propp:d_G_0_G_k}
		has $\vfi\cdot\Gra(W,E,V)=\Gra(W^k,E^k,V)$ if and only if $d(\vfi)=k$.}{It suffices to check that the condition $\vfi\cdot W\in\Gra(W^k,E^k,V)$ is equivalent to $d(\vfi)=k$. 
		Since $\vfi\cdot W$ belongs to $X$, the condition $\vfi\cdot W\in\Gra(W^k,E^k,V)$ is equivalent to the equality  $$\codim_{\vfi\cdot W}((\vfi\cdot W)\cap W)-\codim_W((\vfi\cdot W)\cap W)=k.$$
		We have $\vfi\cdot W=\vfi((W_*)^*)\cap V$ and $W\subset (W_*)^*$, hence $$(\vfi\cdot W)\cap W\subset(\vfi((W_*)^*)\cap V)\cap (W_*)^*=A_0\cap V=A_0\cap W$$
		where $A_0=(W_*)^*\cap\vfi((W_*)^*)$. Moreover, the opposite inclusion is clear, so $(\vfi\cdot W)\cap W=A_0\cap W$. Therefore, we need to prove that
		$$\codim_{\vfi((W_*)^*)\cap V}(A_0\cap W)-\codim_W(A_0\cap W)=d(\vfi).$$
		
		As was shown in the proof of Proposition~A.\ref{propapp:action_X}, there exist $n\in\Zp_{>0}$ and $\alpha_1,~\ldots,~ \alpha_n\in W_*$ such that $\pi_U(\vfi|_{(W_*)^*}(\omega))=\sum_{i=1}^n\omega(\alpha_i)e_i$ for all $\omega\in (W_*)^*$, and $A_0=\{\omega\in (W_*)^*\mid\omega(\alpha_i)=0,~1\leq i\leq n\}$. Hence, $A_0\cap W=\{w\in W\mid\alpha_i(w)=0,~1\leq i\leq n\}$ and
		\begin{equation}\codim_{(W_*)^*}A_0=\codim_W(A_0\cap W)=\dim\pi_U(\vfi((W_*)^*)).\label{formula:A_0_cap_W}
		\end{equation}
		
		It remains to check that $\codim_{\vfi((W_*)^*)\cap V}(A_0\cap W)=\dim\pi_U(\vfi^{-1}((W_*)^*)$. This is done essentially by the same argument as above. Indeed, 
		the argument in the proof of Proposition A.\ref{propapp:action_X} now shows that there exist $m\in\Zp_{>0}$ and $\wt\alpha_i\in W_*$, $1\leq i\leq m$, such that $\pi_U(\vfi^{-1}|_{(W_*)^*}(\omega))=\sum_{i=1}^m\omega(\wt\alpha_i)e_i$ for $\omega\in(W_*)^*$. This, together with fact that $\dim\pi_U(\vfi^{-1}((W_*)^*))<\infty$, implies that, given $(\omega,u)\in V_E$, the condition $\vfi^{-1}((\omega,u))\in(W_*)^*$ is equivalent to a finite system of linear equations on finitely many of the coordinates of $(\omega,u)$. 
		Thus,
		\begin{equation*}
			\begin{split}
				\codim_{\vfi((W_*)^*)\cap V}(A_0\cap W)&=\codim_{\vfi((W_*)^*)}A_0=\codim_{\vfi((W_*)^*)}(\vfi((W_*)^*)\cap (W_*)^*)\\
				&=\codim_{(W_*)^*}((W_*)^*\cap\vfi^{-1}((W_*)^*))=\dim\pi_U(\vfi^{-1}((W_*)^*)).
			\end{split}
		\end{equation*}
		The proof is now complete.}

	\coropappendix{The definition of $d(\vfi)$ does not depend on the choice of $W$ in the following sense: in the definition of $d(\vfi)$ one can replace $W$ by any subspace $W'\subset V$ which is $\widetilde{E}$-commensurable with $W$, and $U$ --- by any direct complement $U'$ of $W'$ in $V$ such that $\codim_{U'} \langle U'\cap \widetilde{E} \rangle_{\mathbb{C}}<\infty$.}{The claim follows directly from Proposition A.\ref{propp:d_G_0_G_k} as the ind-varieties $\Gra(W^k,E^k,V)$ remain unchanged under a replacement $W \leadsto W'$, $U\leadsto U'$.}
	

	Now we are ready to describe the structure of the 
	$ G (V_E, V_{*E})$-action on $X$ in more detail. In particular, we prove that the degree defines a grading on the group $ G (V_E, V_{*E})$. In the next theorem we set $X(k)=\Gra(W^{k},E^{k},V)$ for $k\in \Zp$. Then $X(0)=\Gra(W,E,V)$ and $X=\bigsqcup\limits_{k\in \Zp} X(k)$.
	
	\theopappendix{Given $d\in\Zp$, set \label{theoapp:d_0_G_k_G_k}
		$ G^d(V_E, V_{*E}):=\{\vfi\in  G (V_E, V_{*E})\mid d(\vfi)=d\}$. Then
		\begin{equation*}
			\begin{split}
				&\text{\textup{i)} $\vfi\cdot X(k)=X(k+d(\vfi))$ for all $\vfi\in  G (V_E, V_{*E})$, $k\in\Zp$};\\
				&\text{\textup{ii)} $d(\vfi\circ\vfi')=d(\vfi)+d(\vfi')$ for all $\vfi,~\vfi'\in  G (V_E, V_{*E})$};\\
				&\text{\textup{iii)} $ G^0 (V_E, V_{*E})$ is a normal subgroup of $ G (V_E, V_{*E})$ whose cosets are $ G^d (V_E, V_{*E})$, $d\in\Zp$}.\\
			\end{split}
		\end{equation*}
	}{i) We have $d(\vfi^{-1})=-d(\vfi)$ by the definition of degree. This implies that if $\vfi\cdot X(0)=X(d)$ then $\vfi^{-1}\cdot X(0)=X(-d)$. In particular, since $Sh \cdot X(n)=X(n+1)$ and consequently $Sh^k \cdot X(n)=X(n+k)$, we obtain $Sh^{-k}\cdot X(n)=X(n-k)$ 
		for all $k,~n\in\Zp$. 
		
		We claim that if $\vfi\cdot X(0)= X(0)$ then $\vfi\cdot X(k)= X(k)$ for all $k\in\Zp$. Indeed, assume that $\vfi\cdot X(0)= X(0)$ and pick an integer $k$. We have $(Sh^{-k}\circ\vfi)\cdot X(0)=Sh^{-k}\cdot(\vfi\cdot X(0))=Sh^{-k}\cdot X(0)=X(-k)$. By our observation in the first paragraph, we conclude that $(Sh^{-k}\circ\vfi)^{-1}\cdot X(0)= X(k),$ i.e, $$(\vfi^{-1}\cdot Sh^{-k})\cdot X(0)=\vfi^{-1}\cdot( Sh^k\cdot X(0))=\vfi^{-1}\cdot X(k)= X(k).$$ This clearly implies $\vfi\cdot X(k)= X(k)$.
		
		Now, assume that $\vfi\cdot X(0)=X(d)$, i.e., that $d(\vfi)=d$. Then $( Sh^{-d}\circ\vfi)\cdot X(0)= Sh^{-d}\cdot X(d)= X(0)$. Thus, by the previous paragraph, for all $k\in\Zp$ we have $( Sh^{-d}\circ\vfi)\cdot X(k)= X(k)$, i.e., $$\vfi\cdot X(k)= Sh^d\cdot X(k)=X(k+d).$$
		
		ii) Follows immediately from (i).
		
		iii) Follows immediately from (ii).}
	
	\sstappendix{Matrix realization} 
	Consider the spaces  $V_E=(W_*)^*\oplus U$, $V_{*E}=W_*\oplus U^*$, and assume that the isomorphism $\eta:W_*\rightarrow U$ maps $e^*_{-i}$ to $e_i$ for $i\in \Zp_{>0}$. Then the pairing $V_E\times V_{*E} \rightarrow \Cp$ is identified with the pairing
	\begin{equation*}
		(U^*\oplus U)\times(U^*\oplus U)\to\Cp\colon((\gamma,~y),~(\varkappa,z))\mapsto\langle(\gamma,~y),~(\varkappa,z)\rangle=\gamma(z)+\varkappa(y),~y,~z\in U,~\gamma,~\varkappa\in U^*.
	\end{equation*}
	This allows us to denote the group $ G(V_E,V_{*E})= G(U^*\oplus U,U^*\oplus U)$ simply by $G_U$. 
	
	We now present an explicit matrix realization of the group $G_U$. 
	We have a fixed basis 
	$\{e_1,~e_2,~\ldots\}$ of $U$, 
	and we identify 
	$U^*\oplus U$ with the space $\{(\gamma,~y)\}$ where $\gamma=(\ldots,~\gamma_{-2},~\gamma_{-1})$ are  reverse sequences and $y=(y_{1},~y_2,~\ldots)$  are finitary usual sequences. 
	We consider matrices whose rows and columns are parameterized by the ordered set $\Zp\setminus \{0\}=\Zp_{<0}\sqcup \Zp_{>0}$. 
	Such a matrix $M$ naturally splits into four blocks
	\begin{equation*}
		M=\left(\begin{array}{c|c}
			A & B  \\
			\hline
			C & D  \end{array}\right).
	\end{equation*}
	Let $J$ be the subset of such matrices satisfying the following conditions:
	\begin{equation}
		\begin{split}\label{formula:matrix_group}
			&\text{$\bullet$ each row of $A$ and each column of $D$ is finitary};\\
			&\text{$\bullet$ $C$ is a finitary matrix}.\\
		\end{split}
	\end{equation}
	Matrices from $J$ act on $U^*\oplus U$ via left multiplication: one considers vectors of $U^*\oplus U$ as columns $$(\gamma,y)^{T}=\begin{pmatrix} \gamma^{T} \\ y^T \end{pmatrix},$$ and we have $(\gamma,y)^T \longmapsto \psi(\gamma,y)^{T}$
	for $\psi\in J$.
	
	Given a matrix $\psi\in J$, we let $\overline{\psi}$ denote the reflection of $\psi$ with respect to its antidiagonal.  
	Clearly, $\overline{\overline{\psi}}=\psi$. One 
	checks immediately that $\overline{\psi}$ belongs to $J$, and that $\overline{\psi_1\psi_2}=\overline{\psi_2}~\overline{\psi_1}$. Furthermore, if $\overline{\gamma}:=(\gamma_{-1},~\gamma_{-2},~\gamma_{-3},~\ldots)$ for a reverse sequence $\gamma=(\ldots,~\gamma_{-3},~\gamma_{-2},~\gamma_{-1})$ and $\overline{y}:=(\ldots,~y_3,~y_2,~y_1)$ for a finitary sequence $y=(y_1,~y_2,~y_3,~\ldots)$, then we have $\overline{\psi(\gamma,~y)^T}=\overline{(\gamma,~y)^T}~\overline{\psi}$ for all $(\gamma,~y)\in U^*\oplus U$, where $\overline{(\gamma,~y)^T}=(\overline{y},~\overline{\gamma})$.
	
	Let $\widetilde{G}_U$ be the group of automorphisms of $V_E$ with matrices $\psi$ such that $\psi$, $\psi^{-1}\in J$. 
	
	\theopappendix{The group $G_U$ coincides\label{theo:G_U} with $\wt G_U$.}{First, we prove that $G_U$ is contained in $\wt G_U$. 
		Pick an operator $\vfi\in G_U$. Since $\pi_U\circ (\restr{\vfi}{U})$ is a well-defined operator on the countable-dimensional space $U$
		, it is represented by a (unique) matrix $D=(d_{i,j})_{i,j\in \Zp_{>0}}$ with finitary columns. 
		
		Next, Lemma~A.\ref{lemm:fin_dim_Im} implies that $\dim\Imm(\pi_U\circ(\restr{\vfi}{U^*}))<\infty$. Moreover, 
		according to the proof of Lemma~A.\ref{lemm:fin_dim_Im} there exist $n\in\Zp_{\geq0}$ and vectors $c_1,~\ldots,~c_n\in U$ such that $\pi_U(\restr{\vfi}{U^*}(\gamma))=\sum_{i=1}^n\gamma(c_i)e_i$ for all $\gamma\in U^*$. This means that $\pi_{U}\circ(\restr{\vfi}{U^*})$ can be represented as a finitary matrix $C=(c_{i,-j})_{i,j\in \Zp_{>0}}$, 
		where $c_i=\sum_{j\geq1}c_{i,-j}e_j$.
		
		Let $\pi_{U^*}$ be the projection operator $U^*\oplus U\rightarrow U^*$. The operator $\pi_{U^*}\circ (\restr{\vfi}{U})$ is an operator from the space $U$ of finitary sequences to the space $U^*$ of arbitrary reverse sequences, and is therefore given by a matrix $B=(b_{-i,j})_{i,j\in \Zp_{>0}}$. 
		
		Finally, given $\gamma\in U^*$, we have $\pi_{U^*}(\restr{\vfi}{U^*}(\gamma))(e_i)=\vfi(\gamma)(e_i)=\gamma(\overline{\vfi}(e_i))=\gamma(\pi_U(\overline{\vfi}(e_i)))$. There exist linear functions $\beta_j\in U^*$, $j\in \Zp_{>0}$, such that $\pi_U(\overline{\vfi}(y))=(\beta_1(y),~\beta_2(y),~\ldots)$. Since $\pi_U(\overline{\vfi}(e_i))$ belongs to $U$, there exists $k(i)\in \Zp_{>0}$ such that $\beta_j(e_i)=0$ for $j>k(i)$, i.e., $\pi_U(\overline{\vfi}(e_i))=\sum_{j=1}^{k(i)}\beta_j(e_i)e_j$. Thus, $\pi_{U^*}\circ(\restr{\vfi}{U^*})$ can be represented as a matrix with finitary rows $A=(a_{-i,-j})_{i,j\in \Zp_{>0}}$, where $a_{-i,-j}=\beta_j(e_i)$.
		
		It follows that $\varphi\in 
		G_U$ can be represented by invertible matrices from $L$. This conclusion applies also to the inverse operator $\varphi^{-1}$ by the same argument. 
		
		In order to verify the inclusion $\wt G_U\subset G_U$, let now $\psi$ be a matrix from $L$ such that $\psi^{-1}\in L$. Then the linear operator on $U^*\oplus U$ with matrix $\overline{\psi}$ is the restriction to $U^*\oplus U$ of the operator dual to the operator defined by $\psi$. Indeed,
		\begin{equation}
			\begin{split}
				\langle\overline{\psi}(\gamma,y)^T,(\gamma',y')^T\rangle&=\overline{\overline{\psi}(y,\gamma)^T}(\gamma',y')^T=\overline{(\gamma,y)^T}~\overline{\overline{\psi}}(\gamma',y')^T\\
				&=\overline{(\gamma,y)^T}\psi(\gamma',y')^T=\langle(\gamma,y)^T,\psi(\gamma',y')^T\rangle,
			\end{split}\label{formula:vee_Mackey}
		\end{equation}
		where the three middle terms are products of two or respectively three matrices.
		Moreover, the invertibility of $\overline{\psi}$ follows from the equality $\overline{\psi}^{-1}=\overline{\psi^{-1}}$ which is a consequence of the relation $\overline{\psi_1 \psi_2}=\overline{\psi}_2\overline{\psi}_1$. This shows that $\psi$ determines an operator from $G_U$.
	}
	
	\notaappendix{Finally, the reader will notice that if in the definition of $J$ we replace the ordered set $\Zp \setminus \{0\}=\Zp_{<0}\sqcup \Zp_{>0}$ 
		by any ordered set which is the disjoint union $O_1 \sqcup O_2$ of two linearly ordered countable sets with the condition $r<s$ for $r\in O_1$, $s\in O_2$, and define a group by imposing the conditions (\ref{formula:matrix_group}), we will obtain a group isomorphic to $G_U$.  Of course, in such a setting the matrix of the operator $\overline{\vfi}=\vfi^*|_{U^*\oplus U}$, for a given $\vfi \in   G_U$  with matrix $\psi$, will not be $\overline{\psi}$ as defined above, and its form will depend on the choice of $O_1$ and $O_2$. Moreover, there will also be analogues of the operator $Sh$, and we leave it to the reader to define one.}
	
	The above remark is used in the proof of Theorem \ref{theorem:1} presented in Section \ref{sect:proof_theorem_1_in_the_general_case}.
	
	Finally, using the matrix form of the group $G_U$ given by Theorem~A.\ref{theo:G_U}, it is straightforward to check that $G_U$ is nothing but the group   
	of continuous automorphisms of $U^*\oplus U$ as a Tate space, see, for instance, \cite{Kapranov}. In other words, $G_U$ is the group known as Japanese $GL(\infty)$, and $G^0(U,U^*)$ is simply the connected component of the identity in $G(U,U^*)=G_U$.
	
	\bigskip\textsc{Conflict of interest statement}. The authors declare that they have no conflict of interest.

	\medskip\textsc{Mikhail Ignatyev: Samara National Research University, Ak. Pavlova 1, Samara\break\indent 443011, Russia}
	
	\emph{E-mail address}: \texttt{mihail.ignatev@gmail.com}
	
	\medskip\textsc{Ivan Penkov: Jacobs University Bremen, Campus Ring 1, Bremen 2759, Germany}
	
	\emph{E-mail address}: \texttt{i.penkov@jacobs-university.de}
	
	\newpage

\end{document}